\documentclass[12pt]{amsart}
\usepackage{amsmath}
\usepackage{amssymb}
\usepackage{latexsym}

\newcommand{\Zint}{\mathbb {Z}}    
     
\newcommand{\Rea}{\mathbb {R}}      
\newcommand{\Cplx}{\mathbb {C}}     

\newcommand{\halmos}{\rule{5pt}{5pt}}

\numberwithin{equation}{section}

\newtheorem{prop}{\bf Proposition}[section]
\newtheorem{thm}[prop]{\bf Theorem}
\newtheorem{lemma}[prop]{\bf Lemma}

\newtheorem{cor}[prop]{\bf Corollary}
\newtheorem{conj}{\bf Conjecture}
\newenvironment{remk}{\noindent{\bf Remark}\hskip 5pt}{\hfill{$\Box$}}
\newenvironment{prflemma}{\noindent{\em {\normalsize Proof of the lemma.}}\hskip 5pt}{\hfill{$\Box$}}

\begin{document}
\baselineskip 18pt

\title[Heun equation]
{The Heun equation and the Calogero-Moser-Sutherland system I: the Bethe Ansatz method}
\author{Kouichi Takemura}
\address{Department of Mathematical Sciences, Yokohama City University, 22-2 Seto, Kanazawa-ku, Yokohama 236-0027, Japan.}
\email{takemura@yokohama-cu.ac.jp}


\begin{abstract}
We propose and develop the Bethe Ansatz method for the Heun equation.
As an application, holomorphy of the perturbation for the $BC_1$ Inozemtsev model from the trigonometric model is proved.
\end{abstract}

\maketitle

\section{Introduction}

Olshanetsky and Perelomov proposed a family of integrable quantum systems, which is called the Calogero-Moser-Sutherland system or the Olshanetsky-Perelomov system \cite{OP}.
In the early 90's, Ochiai, Oshima and Sekiguchi classified integrable models of quantum mechanics which are invariant under the action of a Weyl group with certain assumptions \cite{OOS}.
For the case of $B_N$ ($N\geq 3$), the generic model coincides with the $BC_N$ Inozemtsev model.
Oshima and Sekiguchi found that the eigenvalue problem for the Hamiltonian of the $BC_1$ (one particle) Inozemtsev model is transformed to the Heun equation with full parameters \cite{OS}.
Here, the Heun equation is a second-order  Fuchsian differential equation with four regular singular points.

In this paper we will investigate solutions to the Heun equation motivated by the analysis of the Calogero-Moser-Sutherland systems.

In the papers \cite{Tak,KT}, holomorphy of the perturbation for the Calogero-Moser-Sutherland model of type $A_N$ from the trigonometric model was proved, and
in \cite{Tak} some results were obtained by applying the Bethe Ansatz method. Note that the Bethe Ansatz for the model of type $A_N$ was established by Felder and Varchenko \cite{FVthr} by investigating asymptotic behavior of an integral expression of a solution to the KZB equation.

For the case of type $A_1$, the Bethe Ansatz method was found in the textbook of Whittaker and Watson \cite{WW}, although it was covered in the chapter entitled ``Lam\'e's equation'', and the phrase ``Bethe Ansatz'' was not used.
Note that the Lam\'e equation is a special case of the Heun equation \cite{Ron}.

In this paper we will propose and develop the Bethe Ansatz method for the Heun equation, and holomorphy of the perturbation for the $BC_1$ Inozemtsev model from the trigonometric model will be obtained.

The method ``Bethe Ansatz'' appears frequently in physics.
In particular the Bethe Ansatz method has great merits for investigating solvable lattice models \cite{Bax}.
In this paper, we use ``Bethe Ansatz'' in a somewhat restricted sense, similar to the usage in \cite{FVthr,Tak}.
The Bethe Ansatz method replaces the problem of finding eigenstates and eigenvalues of the Hamiltonian (for example (\ref{InoEF})) with a problem of solving transcendental equations for a finite number of variables which are called the Bethe Ansatz equations (for example (\ref{BAeq})).
In our case, we use elliptic functions to describe the Bethe Ansatz equation.

It would be very difficult to solve the transcendental equations explicitly. 
Instead, by applying the trigonometric limit the transcendental equations tend to algebraic equations, which may be more easily solvable.

On the other hand, the Hamiltonian of the $BC_1$ Calogero-Sutherland model appears after applying the trigonometric limit, and the eigenstates for the $BC_1$ Calogero-Sutherland model are described by Jacobi polynomials.
Thus the Jacobi polynomials and the solutions to the Bethe Ansatz equation for the trigonometric case are closely connected, and they are related to the perturbation.
In fact holomorphy of the perturbation for the $BC_1$ Inozemtsev model from the trigonometric model follows from holomorphy of the solutions to the Bethe Ansatz equation with respect to a parameter of a period of elliptic functions.

Now we examine the perturbation. Based on the Hamiltonian $H_{CS}$ of the $BC_1$ Calogero-Sutherland model, the Hamiltonian $H$ of the $BC_1$ Inozemtsev model can be regarded as a perturbed one, i.e. $H=H_{CS}+V_0+\sum _{i=1}^{\infty } V_i(x)p^i$, where $V_0$ is a constant and $V_i(x)$ $(i \in \Zint _{\geq 1})$ are functions.
We can calculate eigenvalues and eigenstates of $H$ as formal power series in $p$ by the perturbation.
For our case, holomorphy of the perturbation can be shown and it implies convergence of the formal power series in $p$ if $|p|$ is sufficiently small.

This paper is organized as follows.
In section \ref{HI}, we clarify the relationship between the Heun equation and the $BC_1$ Inozemtsev system \cite{OS}. 
In section \ref{sec:BAH}, we introduce the Bethe Ansatz method for the $BC_1$ Inozemtsev system, or equivalently for the Heun equation.
In section \ref{sec:invsp}, invariant subspaces of doubly periodic functions are introduced. 
In section \ref{sec:monod}, we obtain an integral representation of eigenfunctions of the $BC_1$ Inozemtsev system.
In section \ref{sec:BAEsigma}, we introduce the Bethe Ansatz method for the $BC_1$ Inozemtsev system and describe the Bethe Ansatz equation.
In section \ref{sec:BAtheta}, we rewrite the Bethe Ansatz equation in terms of theta functions in order to consider the trigonometric limit.
In section \ref{sec:trig}, we consider the trigonometric limit and solve the trigonometric Bethe Ansatz equation.
In section \ref{sec:pert}, relationships between the Bethe Ansatz method and the spectral problem on $L^2$-space for the $BC_1$ Inozemtsev system are clarified. Then holomorphy of the perturbation for the $BC_1$ Inozemtsev model from the trigonometric model is obtained. A key observation is that holomorphy of the perturbation follows from holomorphy of the solutions to the Bethe Ansatz equation.
In section \ref{sec:com}, we give comments, and in the appendix, we note definitions and formulae of elliptic functions.

\section{The Heun equation and the Hamiltonian of the $BC_1$ Inozemtsev system} \label{HI}

The expression of the Heun equation in terms of elliptic functions was already known by Darboux in the 19th century, and this appeared in connection with the soliton theory by Verdier and Treibich \cite{TV}.
In \cite{OS}, the relationship between the eigenfunctions for the $BC_1$ Inozemtsev model and the solutions to the Heun equation was demonstrated.
In this section, we recall the relationship between the $BC_1$ Inozemtsev model and the Heun equation.

The Hamiltonian of the $BC_1$ Inozemtsev model is given as follows:
\begin{equation}
H:= -\frac{d^2}{dx^2} + \sum_{i=0}^3 l_i(l_i+1)\wp (x+\omega_i),
\label{Ino}
\end{equation}
where $\wp (x)$ is the Weierstrass $\wp$-function with periods $(2\omega_1, 2\omega_3)$, $\omega _0=0, \; \omega_2=-(\omega_1 +\omega_3)$ and $l_i$ $(i=0,1,2,3)$ are coupling constants.
In this paper, we will consider eigenvalues and eigenfunctions of the Hamiltonian $H$.
Let $f(x)$ be an eigenfunction of $H$ with an eigenvalue $E$, i.e.  
\begin{equation}
(H-E) f(x)= \left( -\frac{d^2}{dx^2} + \sum_{i=0}^3 l_i(l_i+1)\wp (x+\omega_i)-E\right) f(x)=0.
\label{InoEF}
\end{equation}

We will make a correspondence between equation (\ref{InoEF}) and the Heun equation as was done in \cite{OS}.
Set $z= \wp(x)$ and $e_i=\wp(\omega_i)$ $(i=1,2,3)$.
By a straightforward calculation, we have
$$
\frac{d^2}{dx^2}=\wp '(x)^2\left( \frac{d^2}{dz^2} +\frac{1}{2}\left( \frac{1}{z-e_1}+\frac{1}{z-e_2}+\frac{1}{z-e_3}\right) \frac{d}{dz}\right).
$$
Hence equation (\ref{InoEF}) is equivalent to the equation
\begin{align}
& \left\{ \frac{d^2}{dz^2}+\frac{1}{2}\left( \frac{1}{z-e_1}+\frac{1}{z-e_2}+\frac{1}{z-e_3}\right) \frac{d}{dz} -\frac{1}{4(z-e_1)(z-e_2)(z-e_3)}\cdot \right.
\label{Heun} \\
&  \left. \cdot \left( \tilde{C}+l_0(l_0+1)z+\sum_{i=1}^3 l_i(l_i+1)\frac{(e_i-e_{i'})(e_i-e_{i''})}{z-e_i}\right)\right\}\tilde{f}(z) =0,
\nonumber
\end{align}
where $i', i'' \in \{1,2,3\}$, $i'< i''$, $i'$ and $i''$ are different from $i$, $\tilde{f}(\wp(x))=f(x)$ and $\tilde{C}=-E+\sum_{i=1}^3 l_i(l_i+1)e_i$. Note that $e_1+e_2+e_3=0$.

Equation (\ref{Heun}) has four singular points $e_1,\; e_2, \; e_3$ and $ \infty$ on the Riemann sphere and all singular points are regular. 
The exponents at $z=e_i$ $(i=1,2,3)$ are $(l_i+1)/2$ and $-l_i/2$, and the exponents at $\infty$ are $(l_0+1)/2$ and $-l_0/2$. 

The Heun equation is a standard form of a Fuchsian differential equation with four regular singularities. It is written as 
\begin{equation}
\left( \! \left(\frac{d}{dw}\right) ^2 \! + \left( \frac{\gamma}{w}+\frac{\delta}{w-1}+\frac{\epsilon}{w-t}\right) \frac{d}{dw} +\frac{\alpha \beta w -q}{w(w-1)(w-t)} \right)g(w)=0,
\label{eq:Heun}
\end{equation}
with the condition $\alpha +\beta + 1= \gamma +\delta +\epsilon$.
It is easy to transform an arbitrary Fuchsian differential equation with four regular singularities into the Heun equation. Thereby equation (\ref{Heun}) is transformed to the Heun equation.

Conversely, if a Fuchsian differential equation with four regular singularities is given, we can transform it into equation (\ref{Heun}) with suitable $e_i$ (i=1,2,3) with the condition $e_1+e_2+e_3=0$, $l_i$ $(i=0,1,2,3)$, and $\tilde{C}$ by changing the variable $z \rightarrow \frac{az+b}{cz+d}$ and making the transformation $f \mapsto (z-e_1)^{\alpha_1}(z-e_2)^{\alpha_2}(z-e_3)^{\alpha_3}f$. Thus the relationship between the $BC_1$ Inozemtsev model and the Heun equation is made.

In this paper, we discover properties of the Heun equation by using an expression of elliptic functions especially for the case $l_i \in \Zint_{\geq 0}$ $(i=0,1,2,3)$.
Note that a Fuchsian differential equation with four regular singularities with exponents $(\gamma_i, \delta_i)$ $(i=0,1,2,3)$ satisfying $\gamma_i-\delta_i \in \Zint+\frac{1}{2}$ is transformed to equation (\ref{Heun}) with the condition $l_i \in \Zint_{\geq 0}$.

\section{Bethe Ansatz equation} \label{sec:BAH}
\subsection{The invariant subspace of doubly periodic functions} \label{sec:invsp}
In this subsection we will find finite-dimensional invariant spaces of doubly periodic functions with respect to the Hamiltonian $H$ (see (\ref{Ino})) for the case $l_i \in \Zint_{\geq 0}$ $(i=0,1,2,3)$ and analyze them. In the terminology of the Heun-type equation (\ref{Heun}), doubly periodic functions in $x$  correspond to algebraic functions in the parameter $z(=\wp(x))$.

Let ${\mathcal F}$ be the space spanned by meromorphic doubly periodic functions up to signs, namely
\begin{align}
& {\mathcal F}=\bigoplus _{\epsilon _1 , \epsilon _3 =\pm 1 } {\mathcal F} _{\epsilon _1 , \epsilon _3 }, \label{spaceF} \\
& {\mathcal F} _{\epsilon _1 , \epsilon _3 }=\{ f(x) \mbox{: meromorphic }| f(x+2\omega_1)= \epsilon _1 f(x), \; f(x+2\omega_3)= \epsilon _3 f(x) \} ,
\end{align}
where $(2\omega_1, 2\omega_3)$ are basic periods of elliptic functions.
 
Let $k_i$ $(i=0,1,2,3)$ be the rearrangement of $l_i$ $(i=0,1,2,3)$ such that $k_0\geq k_1 \geq k_2 \geq k_3 (\geq 0)$.
The dimension of the maximum finite-dimensional subspace in ${\mathcal F}$ which is invariant under the action of the Hamiltonian is given as follows.
\begin{thm} \label{periodic}

(i) If the number $k_0+ k_1+ k_2+ k_3 (=l_0+l_1+l_2+l_3)$ is even and $k_0+k_3\geq k_1+k_2$, then the dimension of the maximal finite-dimensional invariant subspace in ${\mathcal F}$ with respect to the action of the Hamiltonian $H$ (see (\ref{Ino})) is $2k_0+1$.

(ii) If the number $k_0+ k_1+ k_2+ k_3$ is even and $k_0+k_3< k_1+k_2$, then the dimension is $k_0+k_1+k_2-k_3+1$.

(iii) If the number $k_0+ k_1+ k_2+ k_3$ is odd and $k_0\geq k_1+k_2+k_3+1 $, then the dimension is $2k_0+1$.

(iv) If the number $k_0+ k_1+ k_2+ k_3$ is odd and $k_0< k_1+k_2+k_3+1$, then the dimension is $k_0+k_1+k_2+k_3+2$.
\end{thm}
\begin{proof}
Firstly, by assuming the existence of the finite-dimensional invariant space $W$, we derive the properties of the space, and later we show its existence.

Let $W$ be a finite-dimensional invariant space in ${\mathcal F}$ and $T^{(i)}$ $(i=0,1,2,3)$ be operators on ${\mathcal F}$ defined by $
T^{(i)}f(x)=f(-x+2\omega_i)$.
From the definitions of the operator $H$ and the periodicity of the function which belongs to ${\mathcal F}$, it follows that $(T^{(i)})^2=1$, $T^{(i)}T^{(j)}=T^{(j)}T^{(i)}$, $T^{(i)}H=HT^{(i)}$ and $T^{(3)}T^{(2)}T^{(1)}T^{(0)}=1$ on ${\mathcal F}$.
Set 
\begin{equation}
V=\sum_{a_0,a_1,a_2,a_3 \in \{0,1\} }(T^{(3)})^{a_3}(T^{(2)})^{a_2}(T^{(1)})^{a_1}(T^{(0)})^{a_0} \cdot W,
\end{equation}
then the space $V$ is finite-dimensional and invariant under the actions of $H$ and $T^{(i)}$ $(i=0,1,2,3)$. From the commutativity of the operators $T^{(i)}$ $(i=0,1,2,3)$, we have the decomposition
\begin{equation}
V= \bigoplus _{\epsilon _0,  \epsilon _1 , \epsilon _2,  \epsilon _3 =\pm 1 \atop{\epsilon _0\epsilon _1\epsilon _2\epsilon _3 =1}}
V_{\epsilon _0, \epsilon _1 , \epsilon _2,  \epsilon _3},
\label{Vdecom0}
\end{equation}
where $V_{\epsilon _0,  \epsilon _1 , \epsilon _2,  \epsilon _3}= \{f(x) \in V |T^{(i)}f(x)= \epsilon _i f(x) \mbox{ for }i=0,1,2,3\} $ $(\epsilon _0,  \epsilon _1 , \epsilon _2,  \epsilon _3 =\pm 1 )$.

Let $f(x)$ be an element of $V_{\epsilon _0,  \epsilon _1 , \epsilon _2,  \epsilon _3}$. From the definition of $T^{(i)}$ and $V_{\epsilon _0,  \epsilon _1 , \epsilon _2,  \epsilon _3}$, the function $f(x+\omega _i)$ $(i=0,1,2,3)$ is odd or even. Hence
\begin{equation} 
f(x)=\sum_{n=M_i}^{\infty } c^{(i)}_{2n+\frac{1-\epsilon_i}{2}}(x-\omega_i)^{2n+\frac{1-\epsilon_i}{2}},
\label{localexp}
\end{equation}
for some $M_i \in \Zint$.
The exponents of the operator $H$ at $x=\omega_i$ are $-l_i$ and $l_i+1$. Let $\alpha _i$ be one of the exponents at $x=\omega_i$ such that $\alpha _i+\frac{1-\epsilon _i}{2}$ is divisible by $2$.

If $c^{(i)}_{2n+\frac{1-\epsilon_i}{2}} \neq 0$ for some $n$ such that $2n+\frac{1-\epsilon_i}{2}<\alpha _i$ then the functions $H^j f(x-\omega_i)$ $(j=0,1,2,\dots )$ are linearly independent and there is a contradiction with finite-dimensionality. Therefore if $2n+\frac{1-\epsilon_i}{2}<\alpha _i$ then $c^{(i)}_{2n+\frac{1-\epsilon_i}{2}}= 0$.

Conversely assume $c^{(i)}_{2n+\frac{1-\epsilon_i}{2}} = 0$ for all $i(\in \{0,1,2,3 \})$ and $n$ such that  $2n+\frac{1-\epsilon_i}{2}<\alpha _i$. Then the function $Hf(x)$ is expressed as $Hf(x)=\sum_{n\in \Zint} \tilde{c}^{(i)}_{2n+\frac{1-\epsilon_i}{2}}(x-\omega_i)^{2n+\frac{1-\epsilon_i}{2}}$ $(i=0,1,2,3)$ and $\tilde{c}^{(i)}_{2n+\frac{1-\epsilon_i}{2}}= 0$ for $n$ such that $2n+\frac{1-\epsilon_i}{2} <\alpha _i$, because $\alpha _i$ is an exponent of the operator $H$ at $x=\omega _i$.

Therefore the space $V_{\epsilon _0,  \epsilon _1 , \epsilon _2,  \epsilon _3}$
consists of functions whose coefficients of the expansion (see (\ref{localexp})) at $x=\omega_i$ $(i\in \{ 0,1,2,3 \} )$ satisfy $c^{(i)}_{2n+\frac{1-\epsilon_i}{2}} = 0$ for all $i$ and $n$ such that $2n+\frac{1-\epsilon_i}{2}<\alpha _i$.
From decomposition (\ref{Vdecom0}) and the condition $\epsilon _0\epsilon _1\epsilon _2\epsilon _3 =1$, the number $\alpha _0+\alpha _1+\alpha _2+\alpha _3$ must be even.

Now we show that the non-zero finite dimensional vector space $V_{\epsilon _0,  \epsilon _1 , \epsilon _2,  \epsilon _3}$ exists iff $\alpha _0+\alpha _1+\alpha _2+\alpha _3 \in 2\Zint_{\leq 0}$.
Let $d=-\frac{\alpha _0 +\alpha _1+\alpha _2+\alpha _3}{2}$ and
\begin{equation}
\tilde{V}_{\epsilon _0,  \epsilon _1 , \epsilon _2,  \epsilon _3} = \left\{
\begin{array}{ll}
\bigoplus _{n=0}^{d} \Cplx \left( \frac{\sigma_1(x)}{\sigma (x)} \right) ^{\alpha_1} \left( \frac{\sigma_2(x)}{\sigma (x)} \right) ^{\alpha_2} \left( \frac{\sigma_3(x)}{\sigma (x)} \right) ^{\alpha_3}\wp(x)^n , & d\in \Zint_{\geq 0} ;\\
\{0\} , & \mbox{otherwise},
\end{array} \right.
\end{equation}
where $\sigma_i(x)$ are the co-sigma functions which are defined in the appendix. Then the space $\tilde{V}_{\epsilon _0,  \epsilon _1 , \epsilon _2,  \epsilon _3}$ is maximum, i.e. $V_{\epsilon _0,  \epsilon _1 , \epsilon _2,  \epsilon _3} \subset \tilde{V}_{\epsilon _0,  \epsilon _1 , \epsilon _2,  \epsilon _3}$ and it is seen that $H \cdot \tilde{V}_{\epsilon _0,  \epsilon _1 , \epsilon _2,  \epsilon _3} \subset \tilde{V}_{\epsilon _0,  \epsilon _1 , \epsilon _2,  \epsilon _3}$. Note that the dimension of the space $\tilde{V}_{\epsilon _0,  \epsilon _1 , \epsilon _2,  \epsilon _3}$ is $d+1$.
Set 
\begin{equation}
\tilde{V}= \bigoplus _{\epsilon _0,  \epsilon _1 , \epsilon _2,  \epsilon _3 =\pm 1 \atop{\epsilon _0\epsilon _1\epsilon _2\epsilon _3 =1}}
\tilde{V}_{\epsilon _0, \epsilon _1 , \epsilon _2,  \epsilon _3},
\label{Vdecom}
\end{equation}
then the space $\tilde{V}$ is the maximum finite-dimensional space spanned by functions in ${\mathcal F}$ and is invariant under the action of the Hamiltonian $H$.
In particular, we have $W\subset V \subset \tilde{V}$.
The dimension of the space $\tilde{V}$ is the summation of the dimension of the spaces $\tilde{V}_{\epsilon _0,  \epsilon _1 , \epsilon _2,  \epsilon _3}$ for $\epsilon _0\epsilon _1 \epsilon _2  \epsilon _3 =1$.

If $l_0+ l_1+ l_2+ l_3$ is even, then the dimension of the space spanned by all finite dimensional invariant spaces in ${\mathcal F}$ is
\begin{align}
&  \left| \frac{l_0+ l_1+ l_2+ l_3}{2}+1 \right| +
\left| \frac{l_0+ l_1- l_2- l_3}{2} \right|  \\
& \; \; \; +\left| \frac{l_0- l_1- l_2+ l_3}{2} \right| +
\left| \frac{l_0- l_1+ l_2- l_3}{2} \right| . \nonumber
\end{align}
If $l_0+ l_1+ l_2+ l_3$ is odd, then the dimension of the space spanned by all finite dimensional invariant spaces in ${\mathcal F}$ is
\begin{align}
&  \left| \frac{l_0+ l_1+ l_2 - l_3+1}{2}\right| +
 \left| \frac{l_0+ l_1- l_2 + l_3+1}{2}\right| \\
& \; \; \; +
 \left| \frac{l_0- l_1+ l_2 + l_3+1}{2}\right| +
 \left| \frac{-l_0+ l_1+ l_2 + l_3+1}{2}\right| .\nonumber
\end{align}
Hence the theorem is obtained.
\end{proof}
\begin{remk}
Theorem \ref{periodic} and its proof is a generalization of \cite[\S 23.41]{WW}.
\end{remk}

We will consider the multiplicities of eigenvalues of the Hamiltonian $H$ on the space spanned by finite dimensional invariant spaces  in ${\mathcal F}$. We will obtain the following theorem, which will be used later.

\begin{thm} \label{thm:dist}
Let $\tilde{V}$ be the space which has appeared in Theorem \ref{periodic} (i.e. (\ref{Vdecom})).
If two of $l_i$ $(i=0,1,2,3)$ are zero, then the roots of the characteristic polynomial of the operator H (see (\ref{Ino})) on the space $\tilde{V}$ are distinct for generic $\omega_1$ and $\omega_3$.
\end{thm}
\begin{proof}
By shifting the variable $x \rightarrow x+\omega_i$ and rearranging the periods, it is sufficient to show the case $l_2=l_3=0$ and $l_0\geq l_1$.
Set $\Phi(z)=(z-e_1)^{\tilde{\alpha}_1}(z-e_2)^{\tilde{\alpha}_2}(z-e_3)^{\tilde{\alpha}_3}$, where $\tilde{\alpha}_i =-l_i/2$ or $(l_i+1)/2$ for each $i(\in \{1,2,3\})$, and $\tilde{f}(z)=\Phi (z)f(z)$. If the function $\tilde{f}(z)$ satisfies equation (\ref{Heun}), then
\begin{align}
& \frac{d^2f(z)}{dz^2}+ \sum_{i=1}^3\frac{2\tilde{\alpha}_i+\frac{1}{2}}{z-e_i} \frac{df(z)}{dz}+\left( \frac{(\tilde{\alpha}_1+\tilde{\alpha}_2+\tilde{\alpha}_3-\frac{l_0}{2})(\tilde{\alpha}_1+\tilde{\alpha}_2+\tilde{\alpha}_3+\frac{l_0+1}{2})z}{(z-e_1)(z-e_2)(z-e_3)} \right.
\label{Heun2} \\
& \left. +\frac{\frac{E}{4}-e_1(\tilde{\alpha}_2+\tilde{\alpha}_3)^2-e_2(\tilde{\alpha}_1+\tilde{\alpha}_3)^2-e_3(\tilde{\alpha}_1+\tilde{\alpha}_2)^2}{(z-e_1)(z-e_2)(z-e_3)}\right)f(z)=0.
\nonumber
\end{align}
Note that the exponents of equation (\ref{Heun2}) at $z=e_2$ are 
$\left\{ \begin{array}{ll}
(0,1/2) & \; (\tilde{\alpha}_2=0) \\
(-1/2,0) & \; (\tilde{\alpha}_2=1/2)
\end{array} \right.$.
Write $f(z)=\sum_{r=0}^{\infty}c_r(z-e_2)^r$, then the recursive relations for coefficients $c_r$ are given as follows:
\begin{align}
& (e_2-e_1)(e_2-e_3)(2\tilde{\alpha}_2+\frac{1}{2})c_{1} +qc_0=0, \label{rec:invsp1} \\
& (r-1+\gamma_1)(r-1+\gamma_2)c_{r-1}+(e_2-e_1)(e_2-e_3)(r+1)(r+2\tilde{\alpha}_2+\frac{1}{2})c_{r+1} \label{rec:invsp} \\
& +(((2\tilde{\alpha}_2+2\tilde{\alpha}_3+r)(e_2-e_1)+(2\tilde{\alpha}_2+2\tilde{\alpha}_1+r)(e_2-e_3))r+q)c_r=0, \nonumber
\end{align}
where $r\geq 1$, $\gamma_1=\tilde{\alpha}_1+\tilde{\alpha}_2+\tilde{\alpha}_3-\frac{l_0}{2}$, $\gamma_2=\tilde{\alpha}_1+\tilde{\alpha}_2+\tilde{\alpha}_3+\frac{l_0+1}{2}$ and $q=\frac{E}{4}-e_1(\tilde{\alpha}_2+\tilde{\alpha}_3)^2-e_2(\tilde{\alpha}_1+\tilde{\alpha}_3)^2-e_3(\tilde{\alpha}_1+\tilde{\alpha}_2)^2+e_2\gamma_1\gamma_2$.
Either $\gamma_1$ or $\gamma_2$ is an integer. Let this integer be denoted by $\tilde{\gamma}$.

The coefficients $c_r$ are polynomial in $E$ of degree $r$. We denote $c_r$ by $c_r(E)$. If $E$ is a solution to the equation $c_{r}(E) =0$ for the positive integer $r=1-\tilde{\gamma}$, then it is seen that $c_{2-\tilde{\gamma}}(E)=0$ from recursive relation (\ref{rec:invsp}), and also $c_{r}(E)=0$ for $r \geq 2-\tilde{\gamma}$.
For the value $E$, the power series $f(z)=\sum_{r} c_r(E) (z-e_2)^r$ is a finite summation and the function $f(\wp(x))$ is meromorphic and doubly periodic. 

If $l_0+l_1$ is even, it occurs that $(\tilde{\alpha}_1, \tilde{\alpha}_2,\tilde{\alpha}_3)=(-l_1/2,0,0)$, $(-l_1/2,1/2,1/2)$, $((l_1+1)/2,1/2,0)$, $ ((l_1+1)/2,0,1/2)$ and the degrees of the polynomials $c_{1-\tilde{\gamma }} (E)$ are $(l_0+l_1+1)/2$, $(l_0+l_1)/2$, $(l_0-l_1)/2$, $(l_0-l_1)/2$ respectively. We denote these polynomials by $c^{(1)}(E),$  $c^{(2)}(E),$  $c^{(3)}(E),$ and $c^{(4)}(E)$. 

Since $\tilde{\alpha_i}$ $(i=0,1,2,3)$ belongs to $\frac{1}{2}\Zint$, the function $\Phi(\wp(x))f(\wp(x))$ is doubly periodic up to signs.
Let $E_1$ (resp. $E_2$) be a solution to $c^{(i_1)}(E)=0$ (resp. $c^{(i_2)}(E)=0$) for some $i_1$ (resp. $i_2$) and let $f_1(z)$ (resp. $f_2(z)$) be the corresponding solution to equation (\ref{Heun2}) which is determined by recursive relations (\ref{rec:invsp1}, \ref{rec:invsp}).
The periodicity of the functions $\Phi(\wp(x))f_1(\wp(x)) $ and  $\Phi(\wp(x))f_2(\wp(x)) $ , (i.e. the signs of $\Phi(\wp(x+\omega_i))f_1(\wp(x+\omega_i))/\Phi(\wp(x))f_1(\wp(x)) $ and $\Phi(\wp(x+\omega_i))f_2(\wp(x+\omega_i))/\Phi(\wp(x))f_2(\wp(x)) $ $(i=1,3)$) are different if $i_1 \neq i_2$.

Assume that $E$ is a common solution to $c^{(i)}(E)=0$ and $c^{(j)}(E)=0$ for $i\neq j$. Let $\tilde{f}^{(i)}(z)$ (resp. $\tilde{f}^{(j)}(z)$) be the solution to differential equation (\ref{Heun2}) related to the value $E$. Then the functions $\tilde{f}^{(i)}(\wp(x))$ and $\tilde{f}^{(j)}(\wp(x))$ are the basis of solutions to (\ref{InoEF}). However, this contradicts the periodicity of the solutions when considering exponents.
Hence, we proved that the equations $c^{(i)}(E)=0$ and $c^{(j)}(E)=0$ $(i\neq j)$ do not have common solutions.

The degree of the polynomial $c^{(1)}(E)c^{(2)}(E)c^{(3)}(E)c^{(4)}(E)$ is equal to the dimension of the space $\tilde{V}$ (see (\ref{Vdecom})).
If we show that the roots of the equation $c^{(i)}(E)=0$ $(i=1,2,3,4)$ are distinct for generic $\omega _1$ and $\omega _3$,
the theorem is proved.

It is enough to show that if $e_1,e_2, e_3$ are real and $(e_2-e_1)(e_2-e_3)<0$ then the roots of the equation $c^{(i)}(E)=0$ $(i=1,2,3,4)$ are real and distinct for the case $(\tilde{\alpha}_1, \tilde{\alpha}_2,\tilde{\alpha}_3)=(-l_1/2,0,0)$, $(-l_1/2,1/2,1/2)$, $((l_1+1)/2,1/2,0)$, $ ((l_1+1)/2,0,1/2)$.

We will show that the polynomial $c_r(E)$ $(1\leq r\leq 1-\tilde{\gamma})$ has $r$ real distinct roots $s_i^{(r)}$ $(i=1,\dots, r)$ such that $s_1^{(r)}<s_1^{(r-1)}<s_2^{(r)}<s_2^{(r-1)}<\dots <s_{r-1}^{(r)}<s_{r-1}^{(r-1)}<s_r^{(r)}$ by induction on $r$. 
The case $r=1$ is trivial.
Let $k \in \Zint _{\geq 1}$ and assume that the statement is true for $r\leq k$.
From the assumption of the induction, $s_1^{(k)}<s_1^{(k-1)}<s_2^{(k)}<s_2^{(k-1)}<\dots <s_{k-1}^{(k-1)}<s_k^{(k)}$.
It is immediate from (\ref{rec:invsp}) that the leading term of the polynomial $c_r(E)$ in $E$ is positive.
 Since $s_i^{(k-1)}$ and $s_{i+1}^{(k-1)}$ satisfy the equation $c_{k-1}(s_i^{(k-1)})=c_{k-1}(s_{i+1}^{(k-1)})=0$ and $s_i^{(k-1)} <s_i^{(k)}< s_{i+1}^{(k-1)}$, the sign of the value $c_{k-1}(s_i^{(k)})$ is $(-1)^{k-i}$. By recursion relation (\ref{rec:invsp}) and the equation $c_{k}(s_i^{(k)})=0$, the sign of the value $c_{k+1}(s_i^{(k)})$ is opposite to that of  $c_{k-1}(s_i^{(k)})$. Combined with the asymptotics of $c_{k+1}(E)$ for $E \rightarrow +\infty$ and $E \rightarrow -\infty$ and $\deg_E c_{k+1}(E)=k+1$, it follows that the polynomial $c_{k+1}(E)$  has $k+1$ real distinct roots and the inequality  $s_1^{(k+1)}<s_1^{(k)}<s_2^{(k+1)}<s_2^{(k)}<\dots <s_{k}^{(k)}<s_{k+1}^{(k+1)}$ is satisfied.

For the case that $l_0+l_1$ is odd, the proof is similar. 
\end{proof}
The following proposition was essentially proved in the proof of the previous theorem.
\begin{prop}
Suppose $l_2=l_3=0$ and $l_0, l_1\in \Zint_{\geq 0}$. If $e_1,e_2, e_3$ are real and $(e_2 -e_1)(e_3-e_1) >0$, then the roots of the characteristic polynomial of operator $H$ (\ref{Ino}) on the space $\tilde{V}$ are real and distinct.
\end{prop}
\begin{proof}
It is sufficient to prove for the case $l_0\geq l_1$.
If $e_1,e_2, e_3$ are real then the condition $(e_2 -e_1)(e_3-e_1) >0$ is equivalent to ($(e_2 -e_1)(e_2-e_3) <0$ or $(e_3 -e_1)(e_3-e_2) <0$).
For the case $(e_2 -e_1)(e_2-e_3) <0$, the proposition is proved in the proof of the previous theorem. For the case $(e_3 -e_1)(e_3-e_2) <0$, the proof is similar.
\end{proof}

\subsection{Results from the method of monodromy} \label{sec:monod}
In this subsection, we present an integral expression of solutions to equation (\ref{InoEF}) and related results for the case $l_i \in \Zint_{\geq 0}$ $(i=0,1,2,3)$.
First, we show a proposition which is related to the monodromy of solutions to equation (\ref{InoEF}) for the case $l_i \in \Zint_{\geq 0}$ $(i=0,1,2,3)$.
\begin{prop} \label{prop:locmonod}
If $l_0, l_1, l_2, l_3  \in \Zint _{\geq 0}$ then the monodromy matrix of equation (\ref{InoEF}) around each singular point is a unit matrix.
\end{prop} 
\begin{proof}
Each solution to equation (\ref{InoEF}) has singular points only at $x=n_1 \omega_1 + n_3 \omega_3$ $(n_1, n_3 \in \Zint)$.
Due to periodicity, it is sufficient to consider the case $a=\omega _i$ $(i=0,1,2,3)$.
The exponents at the singular point $a=\omega_i$ are $-l_i$ and $l_{i}+1$, hence there exist solutions of the form $f_{a,1}(x)=(x-a)^{-l_i}(1+\sum_{k=1}^{\infty}a_k x^{2k})$ and $f_{a,2}(x)=(x-a)^{l_i+1}(1+\sum_{k=1}^{\infty}a'_k x^{2k})$, because equation (\ref{InoEF}) lacks a first-order differential term and because the functions $\wp (x+\omega_j )$ $(j=0,1,2,3)$ are even and because the numbers $l_j+1-(-l_j)$ are odd.
Then it is obvious that the monodromy matrix for the functions $f_{a,1}(x)$ and $f_{a,2}(x)$ around the point $x=a$ is a unit matrix, hence the local monodromy matrix of equation (\ref{InoEF}) around each singular point is also a unit matrix.
\end{proof}

We explain how to choose particular solutions to equation (\ref{InoEF}) by applying Proposition \ref{prop:locmonod}.

Let $M_i$ $(i=1,3)$ be the transformations of solutions to differential equation (\ref{InoEF}) obtained by the analytic continuation $x \rightarrow x+2\omega _i$. Since all the local monodromy matrices are unit, the transformations $M_i$ do not depend on the choice of paths. From the fact that the fundamental group of the torus is commutative, we have $M_1 M_3=M_3 M_1$.
The operators $M_i$ act on the space of solutions to equation (\ref{InoEF}) for each $E$, which is two dimensional.
By the commutativity $M_1 M_3=M_3 M_1$, there exists a joint eigenvector $\Lambda (x,E)$ for the operators $M_1$ and $M_3$. From Proposition \ref{prop:locmonod}, the function $\Lambda (x,E)$ is single-valued and satisfies equations $(H-E)\Lambda (x,E)=0$, $M_1\Lambda (x,E)=m_1\Lambda (x,E)$ and $M_3\Lambda (x,E)=m_3\Lambda (x,E)$ for some $m_1,m_3 \in \Cplx \setminus \{0\}$.
We immediately obtain relations $(H-E)\Lambda (-x,E)=0$, $M_1\Lambda (-x,E)=m_1^{-1}\Lambda (-x,E)$ and $M_3\Lambda (-x,E)=m_3^{-1}\Lambda (-x,E)$.

Now consider the function $\Xi ^* (x,E)=\Lambda (x,E)\Lambda (-x,E)$. This function is even (i.e. $\Xi ^*(x,E)=\Xi ^*(-x,E)$), doubly periodic (i.e. $\Xi ^*(x+2\omega _1,E)=\Xi ^*(x+2\omega _3,E)=\Xi ^*(x,E)$), and satisfies the equation
\begin{align}
& \left( \frac{d^3}{dx^3}-4\left( \sum_{i=0}^3 l_i(l_i+1)\wp (x+\omega_i)-E\right)\frac{d}{dx} \right. \nonumber \\
& \left. -2\left(\sum_{i=0}^3 l_i(l_i+1)\wp '(x+\omega_i)\right) \right) \Xi ^*(x,E)=0,
\nonumber
\end{align}
which products of any pair of solutions to equation (\ref{InoEF}) satisfy.

The following proposition is used to determine an appropriate constant multiplication of the function $\Xi ^* (x,E)$.

\begin{prop} \label{prop:prod}
If $l_0, l_1, l_2, l_3 \in \Zint _{\geq 0}$ then the equation
\begin{align}
& \left( \frac{d^3}{dx^3}-4\left( \sum_{i=0}^3 l_i(l_i+1)\wp (x+\omega_i)-E\right)\frac{d}{dx} \right. \label{prodDE} \\
& \left. -2\left(\sum_{i=0}^3 l_i(l_i+1)\wp '(x+\omega_i)\right) \right) \Xi (x,E)=0,
\nonumber
\end{align}
has a nonzero doubly periodic solution which has the expansion
\begin{equation}
\Xi (x,E)=c_0(E)+\sum_{i=0}^3 \sum_{j=0}^{l_i-1} b^{(i)}_j (E)\wp (x+\omega_i)^{l_i-j},
\label{Fx}
\end{equation}
where the coefficients $c_0(E)$ and $b^{(i)}_j(E)$ are polynomials in $E$, they do not have common divisors, and the polynomial $c_0(E)$ is monic.
Moreover the function $\Xi (x,E)$ is determined uniquely. 
Set $g=\deg_E c_0(E)$, then $\deg _E b^{(i)}_j(E)<g$ for all $i$ and $j$.
\end{prop}
\begin{proof}
Recall that the function $\Xi ^* (x,E)=\Lambda (x,E)\Lambda (-x,E)$ is even and doubly periodic, and it satisfies equation (\ref{prodDE}).

Here we show the following lemma.
\begin{lemma} \label{lem:onedim}
Suppose $l_0, l_1, l_2, l_3 \in \Zint _{\geq 0}$ and fix $E\in \Cplx$.
If equation (\ref{InoEF}) does not have solutions in ${\mathcal F}$, then the dimension of the space of even doubly periodic functions which satisfy equation (\ref{prodDE}) is one.
\end{lemma}
\begin{prflemma}
Since the dimension of even (or odd) functions  (not necessarily doubly periodic) which satisfy equation (\ref{InoEF}) is one and solutions to equation (\ref{prodDE}) are spanned by products of two solutions to (\ref{InoEF}), the dimension of even functions which satisfy equation (\ref{prodDE}) is two.

Suppose that the function $\Lambda (x,E)$ is proportional to $\Lambda (-x,E)$, then by definition one has $m_1=m_1^{-1}$ and $m_3=m_3^{-1}$, that is  $m_1^2=m_3^2=1$, and it contradicts to the assumption of the lemma. Therefore the functions $\Lambda (x,E)$ and $\Lambda (-x,E)$ are linearly independent.

If the dimension of even doubly periodic functions which satisfy equation (\ref{prodDE}) is no less than two, all even solution to (\ref{prodDE}) must be doubly periodic. Hence the function $\Lambda (x,E)^2 +\Lambda (-x,E)^2$ is doubly periodic. Since $\Lambda (x,E) \not \in {\mathcal F}$, it follows that $m_1^2 \neq 1$ or $m_3^2 \neq 1$. For example, suppose $m_1^2\neq 1$. From the periodicity of $ \Lambda (x,E)$, $\Lambda (x,E)^2 +\Lambda (-x,E)^2 = \Lambda (x+1,E)^2 +\Lambda (-x-1,E)^2 = m_1^2 \Lambda (x,E)^2 + m_1^{-2}\Lambda (-x,E)^2$ which contradicts the linear independence of $\Lambda (x,E)$ and $\Lambda (-x,E)$.
Therefore the dimension of even doubly periodic functions which satisfy equation (\ref{prodDE}) is less than two.
Since the function $\Xi ^* (x,E)$ is an even doubly periodic function that satisfies (\ref{prodDE}), we obtain the lemma.
\end{prflemma}

We continue the proof of Proposition \ref{prop:prod}. Since the function $\Xi ^* (x,E)$ is an even doubly periodic function that satisfies the differential equation (\ref{prodDE}) and the exponents of equation (\ref{prodDE}) at $x=\omega _i$ are $-2l_i,1,2l_i+2$ , it admits the following expansion:
\begin{equation}
\Xi ^* (x,E)=\tilde{c}_0(E)+\sum_{i=0}^3 \sum_{j=0}^{l_i-1} \tilde{b}^{(i)}_j (E)\wp (x+\omega_i)^{l_i-j}.
\label{Xitildex}
\end{equation}
Substituting (\ref{Xitildex}) for (\ref{prodDE}), we obtain relations for the coefficients $\tilde{c}_0(E)$ and $\tilde{b}^{(i)}_j (E)$. On the relations, the coefficient $\tilde{b}^{(i)}_j (E)$ is expressed as a polynomial in $\tilde{b}^{(i)}_{j'} (E)$ $(j'<j)$ and $E$, and the coefficient $\tilde{c}_0(E)$ is expressed as a polynomial in $\tilde{b}^{(i)}_{j'} (E)$ $(i=0,1,2,3; \: 0\leq j\leq l_i-1)$ and $E$. 

Through a multiplicative change of scale given by
\begin{equation}
\Xi(x,E)=\frac{c_0(E)}{\tilde{c}_0(E)}\Xi ^*(x,E)=c_0(E)+\sum_{i=0}^3 \sum_{j=0}^{l_i-1} b^{(i)}_j (E)\wp (x+\omega_i)^{l_i-j},
\label{Fx2}
\end{equation}
we can choose the coefficients in such a way that $c_0(E)$ and $b^{(i)}_j(E)$ are polynomials in $E$, they do not have common divisors, and the polynomial $c_0(E)$ is monic.

If there is a function  in ${\mathcal F}$ that satisfies $(H-E)f(x)=0$, we have $f(x) \in \tilde{V}_{\epsilon _0, \epsilon _1 , \epsilon _2,  \epsilon _3}$ (see (\ref{Vdecom})) for some $\epsilon _0, \epsilon _1 , \epsilon _2,  \epsilon _3 \in \{\pm 1 \}$, and the value $E$ is an eigenvalue of $H$ on the space $V_{\epsilon _0, \epsilon _1 , \epsilon _2,  \epsilon _3}$. From Theorem \ref{periodic} and Lemma \ref{lem:onedim}, any even doubly periodic function that satisfies equation (\ref{prodDE}) is determined uniquely up to constant multiplication except for finitely many $E$. Hence the polynomials $c_0(E)$ and $b^{(i)}_j(E)$ are determined uniquely.

Let $\Xi(x,E) = \sum_{k=0}^{g'} a_{g'-k}(x) E^k$, where $g'$ is the maximum of the degrees of $c_0(E)$ and $b^{(i)}_{j'} (E)$ $(i=0,1,2,3; \: 0\leq j\leq l_i-1)$ in $E$. Substituting this into (\ref{prodDE}) and computing the coefficient of $E^{g'+1}$ yields $\frac{d}{dx}a_0(x)=0$, and therefore we have $\deg _E b^{(i)}_j(E)<\deg _E c_0(E)$, $g'=g$, and $a_0(x)=1$. 
\end{proof}

We will get an integral representation of $\Lambda (x,E)$ by using the function $\Xi (x,E)$. The following argument is motivated by the book of Whittaker and Watson \cite[\S 23.7]{WW}. 

Since the functions $\Lambda (x,E)$ and $\Lambda (-x,E)$ satisfy differential equation (\ref{InoEF}), we have $\frac{d}{dx} \left( \Lambda (-x,E) \frac{d}{dx} \Lambda (x,E) -\Lambda (x,E) \frac{d}{dx} \Lambda (-x,E) \right) = \Lambda (-x,E) \frac{d^2}{dx^2} \Lambda (x,E) -  \Lambda (x,E) \frac{d^2}{dx^2}\Lambda (-x,E)= 0$.
Therefore $\Lambda (-x,E) \frac{d}{dx} \Lambda (x,E) -\Lambda (x,E) \frac{d}{dx} \Lambda (-x,E) = 2C(E)$ and $\frac{d \log \Lambda (x,E)}{dx} - \frac{d \log \Lambda (-x,E)}{dx}=\frac{2C(E)}{\Xi (x,E)}$ for some constant $C(E)$.
From the equation $\frac{d \log \Lambda (-x,E)}{dx}+\frac{d \log \Lambda (x,E)}{dx}=\frac{1}{\Xi (x,E)}\frac{d\Xi (x,E)}{dx}$, we have 
\begin{equation}
\frac{d \log \Lambda(\pm x ,E)}{dx}=\frac{1}{2\Xi (x,E)}\frac{d\Xi (x,E)}{dx}\pm  \frac{C(E)}{\Xi (x,E)}.
\label{logfx}
\end{equation}
which is integrated to yield
\begin{equation}
\Lambda (\pm x, E)=\sqrt{\Xi (x,E)}\exp \int \frac{\pm C(E)dx}{\Xi (x,E)}.
\label{integ}
\end{equation}
This formula is the integral representation of the solution to equation (\ref{InoEF}).

It remains to determine the constant $C(E)$. Differentiating equation (\ref{logfx}) yields
\begin{align}
& \frac{1}{\Lambda (x,E)}\frac{d^2\Lambda (x,E)}{dx^2}-\left(\frac{1}{\Lambda (x,E)}\frac{d\Lambda (x,E)}{dx} \right)^2 \label{intdiff} \\
& =\frac{1}{2\Xi (x,E)}\frac{d^2\Xi (x,E)}{dx^2}-\frac{1}{2}\left(\frac{1}{\Xi (x,E)}\frac{d\Xi (x,E)}{dx} \right)^2-\frac{C(E)}{\Xi (x,E)^2}\frac{d\Xi (x,E)}{dx} \nonumber .
\end{align}
Set $Q(E)= -C(E)^2$. By substituting (\ref{InoEF}) and (\ref{logfx}) into (\ref{intdiff}), we obtain the formula
\begin{align}
 Q(E)= -C(E)^2= & \Xi (x,E)^2\left( E- \sum_{i=0}^3 l_i(l_i+1)\wp (x+\omega_i)\right) \label{const} \\
& +\frac{1}{2}\Xi (x,E)\frac{d^2\Xi (x,E)}{dx^2}-\frac{1}{4}\left(\frac{d\Xi (x,E)}{dx} \right)^2. \nonumber
\end{align}
It follows that $Q(E)$ is a monic polynomial in $E$ of degree $2g+1$ from the expression for $\Xi (x,E)$ given by (\ref{Fx}) and Proposition \ref{prop:prod}.
In summary, we have the following proposition.
\begin{prop} 
Let $\Xi (x,E)$ be the doubly periodic function defined in Proposition \ref{prop:prod} and $Q(E)$ be the monic polynomial defined in (\ref{const}).
Then the function 
\begin{equation}
\Lambda (\pm x, E)=\sqrt{\Xi (x,E)}\exp \int \frac{\pm \sqrt{-Q(E)}dx}{\Xi (x,E)}.
\label{integ1}
\end{equation}
is a solution to the differential equation (\ref{InoEF}).
\end{prop}

From the previous proposition, if $E$ is a solution to the equation $Q(E)=0$, then the corresponding function $\Lambda (x,E)$ is an element of ${\mathcal F}$.

Conversely, we suppose $\Lambda (x,E) \in {\mathcal F}$. For this case, the function $\Lambda (-x,E)$ has the same periodicity as $\Lambda (x,E)$.
If the function $\Lambda (x,E)$ is neither even nor odd, then the function $\Lambda (x,E)+\Lambda (-x,E)$ is non-zero, even, and doubly periodic up to signs, and the non-zero function $\Lambda (x,E)- \Lambda (-x,E)$ has the same periodicity.
The function $\Lambda (x,E)+\Lambda (-x,E)$ belongs to the space $\tilde{V}_{\epsilon _0, \epsilon _1 , \epsilon _2,  \epsilon _3}$ (see (\ref{Vdecom})) for some $\epsilon _0, \epsilon _1 , \epsilon _2,  \epsilon _3 \in \{\pm 1 \}$. 
Since the function $\Lambda (x,E)- \Lambda (-x,E)$ has the same periodicity and the opposite parity as $\Lambda (x,E)+\Lambda (-x,E)$, we have $\Lambda (x,E)- \Lambda (-x,E) \in \tilde{V}_{-\epsilon _0, -\epsilon _1 , -\epsilon _2,  -\epsilon _3}$.
Let $\alpha _i\in \{ -l_i, l_i +1 \}$ be the number such that $\alpha _i+\frac{1-\epsilon _i}{2}$ is divisible by $2$,
then we have $\dim \tilde{V}_{\epsilon _0, \epsilon _1 , \epsilon _2,  \epsilon _3} = \max (0, 1-\frac{\alpha _0+\alpha _1+\alpha _2+\alpha _3}{2})$ and $\dim \tilde{V}_{-\epsilon _0, -\epsilon _1 , -\epsilon _2,  -\epsilon _3} = \max (0, \frac{\alpha _0+\alpha _1+\alpha _2+\alpha _3}{2} -1)$. Hence either the space $\tilde{V}_{\epsilon _0, \epsilon _1 , \epsilon _2,  \epsilon _3}$ or the space $\tilde{V}_{-\epsilon _0, -\epsilon _1 , -\epsilon _2,  -\epsilon _3}$ is $\{ 0 \}$, which contradicts the property $\Lambda (x,E) \neq \pm \Lambda (-x,E)$. Therefore the function $\Lambda (x,E)$ must be odd or even.
If the function $\Lambda (x,E)$ is odd or even, the functions $\Lambda (x,E)$ and $\Lambda (-x,E)$ are linearly dependent.
Thus the Wronskian $\Lambda (-x,E) \frac{d}{dx}\Lambda (x,E) -\Lambda (x,E)\frac{d}{dx}\Lambda (-x,E)$ is equal to $0$ and $C(E)=Q(E)=0$.

Therefore the following theorem is proved:
\begin{thm} \label{thm:P(E)}
The equation (\ref{InoEF}) has a solution in the space ${\mathcal F}$ (\ref{spaceF}) if and only if the value $E$ satisfies the equation $Q(E)=0$.
\end{thm}

Let us count the degree of the polynomial $Q(E)$ and discuss the relationship to the maximum finite-dimensional invariant subspace $\tilde{V}$ (see (\ref{Vdecom})) in  ${\mathcal F}$.

\begin{conj} \label{conj1}
(i) The degree of the polynomial $Q(E)$ in $E$ is equal to the dimension of the space $\tilde{V}$. (See (\ref{Vdecom}) and Theorem \ref{periodic})

(ii) For generic $\omega_1 $ and $\omega_3$, roots of the equation $Q(E)=0$ are distinct.
\end{conj}

\begin{prop} \label{propconj1}
Conjecture \ref{conj1} is true if two of $l_i$ $(i=0,1,2,3)$ are zero.
\end{prop}
\begin{proof}
It is sufficient to show the case $l_2=l_3=0$ and $l_0\geq l_1$.

The function $\Xi (x,E)$ admits expression (\ref{Fx}), but we now use the other expression,
\begin{equation}
\Xi (x,E)=\frac{1}{(z-e_1)^{l_1}}\sum_{j=0}^{l_0+l_1}a_j (E)(z-e_1)^j,
\end{equation}
where $z=\wp(x)$.
We are going to obtain a recursive relation for the coefficients $a_j(E)$.

Since $\Xi (x,E)$ satisfies equation (\ref{prodDE}), we find that
\begin{align}
& j(j-2l_0-1)(2j-2l_0-1)a_j (E)\label{recrelF}\\
& +2(j-l_0-1)(-E+(3j^2-6(l_0+1)j+2l_0^2+5l_0-l_1^2-l_1+3)e_1)a_{j-1}(E) \nonumber \\
& +(2j-2l_0-3)(j-2-l_0-l_1)(j-1-l_0+l_1)(e_1-e_2)(e_1-e_3)a_{j-2}(E)=0 ,\nonumber
\end{align}
where $j=1,2,\dots ,l_0+l_1+1$ and $a_{-1}(E)=a_{l_0+l_1+1}(E)=0$.

By recursive relation (\ref{recrelF}), we have\\
(a) All coefficients $a_j(E)$ $(j=1,\dots l_0+l_1)$ are divisible by $a_0(E)$ as a polynomial in $E,e_1,e_2,e_3$. Hence we may set $a_0(E)=1$.\\
(b) The coefficient $a_j(E)$ $(j=0,1,\dots ,l_0+l_1)$ is a polynomial in variables $E, e_1, e_2, e_3$ and the homogeneous degree of $a_j(E)$ w.r.t. $E, e_1, e_2, e_3$ is $j$. \\
(c) $\deg_E a_j(E) \leq j$. \\
(d) If $j\leq l_0$ then $\deg_E a_j(E) = j$.

Put $j=l_0+l_1+1$ in (\ref{recrelF}), then $\deg_Ea_{l_0+l_1-1}(E)\leq \deg_Ea_{l_0+l_1} (E) +1$. By putting $j=l_0+l_1-k+1$ $(k=1,\dots ,l_1-1)$ in (\ref{recrelF}), we also obtain $\deg_Ea_{l_0+l_1-k-1}(E)\leq \deg_Ea_{l_0+l_1-k}(E)+1$.
It follows that $\deg_E a_j(E) \leq l_0$ for all $j$, $\deg_E a_j (E)= l_0 \Leftrightarrow  j=l_0 $, and the coefficient of $E^{l_0}$ in $a_{l_0}(E)$ does not depend on $e_1$, $e_2$, $e_3$.

From relation (\ref{const}), it is seen that $\deg_E Q(E) = 2l_0+1$ and the coefficient of $E^{2l_0+1}$ in $Q(E)$ does not depend on $e_1$, $e_2$, $e_3$.
From Theorem \ref{thm:P(E)}, equation (\ref{InoEF}) has a solution  in ${\mathcal F}$ iff the value $E$ satisfies the equation $Q(E)=0$, and from Theorem \ref{periodic}, $\dim \tilde{V} = 2l_0+1$ (see (\ref{Vdecom})).
From Theorem \ref{thm:dist}, the roots of the characteristic polynomial of the operator $H$ (see (\ref{Ino})) on the space $\tilde{V}$ are distinct for generic $\omega_1$ and $\omega_3$.
Hence the characteristic polynomial of the operator $H$ on the space $\tilde{V}$ coincides with the polynomial $Q(E)$ up to constant multiplication.

Therefore we obtain the proposition.
\end{proof}

The following proposition is proved case by case.
\begin{prop}
Conjecture \ref{conj1} is true for the case $l_0+l_1+l_2+l_3\leq 8$.
\end{prop}

\subsection{Bethe Ansatz equation in terms of sigma functions} \label{sec:BAEsigma}

In the previous subsection, we introduced the elliptic function $\Xi (x,E)$ and the polynomial $Q(E)$.
In this subsection, we will propose the Bethe Ansatz method.
We postulate the Ansatz (assumption) that the eigenfunction of the Hamiltonian (\ref{InoEF}) is written as (\ref{BV2}) with auxiliary variables $(t_1, \dots ,t_l)$.
The Bethe Ansatz method replaces the spectral problem (see (\ref{InoEF})) with solving transcendental equations (\ref{BAeq}) for a finite number of variables $(t_1, \dots ,t_l)$, which we call the Bethe Ansatz equation. 
It will be shown that the basis of solutions to equation (\ref{InoEF}) is written as the form satisfying Ansatz (\ref{BV2}) except for finitely many eigenvalues $E$. The precise statement will be given in Theorem \ref{BAthm}.
Note that the case of the Lam\'e equation (i.e. the case $l_1=l_2=l_3=0$) was treated in \cite{WW}.

Throughout this subsection, we assume $Q(E) \neq 0$. As was discussed in section \ref{sec:monod}, the functions $\Lambda (x,E)$ and $\Lambda (-x,E)$ are linearly independent.

For simplicity, assume $l_0\neq 0$.
The function $\Xi (x,E)$ is an even doubly periodic function which may have poles of degree $2l_i$ at $x=\omega_i$ $(i=0,1,2,3)$.
Therefore we have the following expression;
\begin{equation}
\Xi (x,E)=\frac{a_0^{(0)}\prod_{j=1}^{l}(\wp(x)-\wp(t_j))}{(\wp(x)-e_1)^{l_1}(\wp(x)-e_2)^{l_2}(\wp(x)-e_3)^{l_3}} \label{Fxtj}
\end{equation}
for some values $t_1,\dots ,t_l$, where $l=l_0+l_1+l_2+l_3$.
Note that the values $t_1,\dots ,t_l$ depend on $E$.

\begin{prop} \label{prop:tjdist}
Assume $Q(E) \neq 0$, then the values $t_j$, $-t_j$, $\omega_i$ $(j=1,\dots ,l,\; i=0,1,2,3,$ see $(\ref{Fxtj}))$ are mutually distinct up to periods of elliptic functions.
\end{prop}
\begin{proof}
(a) We show $t_j \neq \omega _i$ up to periods of elliptic functions for each $j \in \{1,\dots ,l \}$ and $i\in \{0,1,2,3 \}$.

Since the function $\Xi (x,E)$ is an even doubly periodic function and satisfies equation (\ref{prodDE}), it has a pole of degree $2l_i$ or a zero of degree $2l_i+2$ at $x=\omega_i$ $(i=0,1,2,3)$.
Suppose that the function $\Xi (x,E)$ has a zero of degree $2l_i+2$ at $x=\omega_i$ for some $i\in \{0,1,2,3 \}$. Then the functions $\Lambda (x,E)$ and $\Lambda (-x,E)$ have zeros of degree $l_i+1$ at $x=\omega_i$, because $\Xi (x,E)$ is proportional to the product $\Lambda(x,E) \Lambda(-x,E)$, $\Lambda(x,E)$ and $ \Lambda(-x,E)$ satisfy differential equation (\ref{InoEF}) and their exponents at $x=\omega_i$ are $-l_i$ and $l_i+1$.
From the condition $Q(E) \neq 0$, the functions $\Lambda (x,E)$ and $\Lambda (-x,E)$ are linearly independent and form a basis of solutions to equation (\ref{InoEF}). It follows that every solution to equation (\ref{InoEF}) has a zero of degree $l_i+1$ at $x=\omega_i$ and this contradicts the fact that one of the exponents at $x=\omega_i$ is $-l_i$.
Therefore the function $\Xi (x,E)$ has a pole of degree $2l_i$ at $x=\omega_i$ $(i=0,1,2,3)$.

If $t_j= \omega_i$ or $- \omega_i$ for some $j$ and $i$, then from expression (\ref{Fxtj}) it contradicts the degree of the pole at $x=\omega_i$. 

(b) Now we show  $t_j \neq t_{j'}$ $(j\neq j', \; j,j' \in \{1,\dots ,l \})$ and $t_j \neq -t_{j'}$ $( j,j' \in \{1,\dots ,l \})$ up to the periods of elliptic functions.
If $t_j=-t_j$, then we have $t_j=\omega_i$ for some $i\in \{ 0,1,2,3\}$ from the condition $2t_j=0$ mod $2\omega_1\Zint +  2\omega_3\Zint$ and it is reduced to the case (a).

Assume ($t_j = t_{j'}$ $(j\neq j', \; j,j' \in \{1,\dots ,l \})$ or $t_j = -t_{j'}$ $( j,j' \in \{1,\dots ,l \})$) and  $t_j \neq \omega _i$ $(j \in \{1,\dots ,l \} ,i\in \{0,1,2,3 \})$, then the function $\Xi (x,E)$ has zeros of degree no less than 2 at two different points $x=t_j$ and $x=-t_j$. 
The degree of zeros of the functions $\Lambda (x,E)$ and $\Lambda (-x,E)$ at points $x=t_j$ and $x=-t_j$ are $0$ or $1$, because $\Lambda (x,E)$ and $\Lambda (-x,E)$ satisfy differential equation (\ref{InoEF}) which is holomorphic at the points $x=\pm t_j$.
Since $\Xi (x,E)$ is proportional to $\Lambda(x,E) \Lambda(-x,E)$, there is a contradiction with the degree of zeros of the function $\Xi (x,E)$.
Therefore the proposition is proved.
\end{proof}

Set $z=\wp(x)$ and $z_j=\wp(t_j)$. From formula (\ref{const}), we obtain
$$
\left( \frac{d\Xi (x,E)}{dz}\right) _{z=z_j}^2=\frac{4C(E)^2}{\wp'(t_j)^2}.
$$
We fix the signs of $t_j$ by taking 
$$
\left( \frac{d\Xi (x,E)}{dz}\right) _{z=z_j}=\frac{2C(E)}{\wp'(t_j)}.
$$
If we put $2C(E)/\Xi (x,E)$ into partial fractions, it is seen that
\begin{align}
& \frac{2C(E)}{\Xi (x,E)}=\sum_{j=1}^l \frac{\frac{2C(E)}{\left( \frac{d\Xi}{dz}\right) _{z=z_j}}}{\wp (x)-\wp (t_j)}  = \sum_{j=1}^l \left(\zeta(x-t_j)-\zeta(x+t_j)+2\zeta(t_j) \right) \nonumber ,
\end{align}
where $\zeta (x)$ is the Weierstrass zeta function. It follows from formula (\ref{integ}) that
\begin{align}
& \Lambda (x,E)=\sqrt{\frac{a_0^{(0)}\prod_{j=1}^{l}(\wp(x)-\wp(t_j))}{(\wp(x)-e_1)^{l_1}(\wp(x)-e_2)^{l_2}(\wp(x)-e_3)^{l_3}}} \label{BA} \\
& \exp\left( \frac{1}{2}\sum_{j=1}^l \left( \log \sigma(t_j+x)-\log \sigma (t_j-x) -2x\zeta(t_j)\right) \right) \nonumber \\
& = \frac{\sqrt{a_0^{(0)}}\prod_{j=1}^l \sigma(x+t_j)}{\sigma(x)^{l_0}\sigma_1(x)^{l_1}\sigma_2(x)^{l_2}\sigma_3(x)^{l_3}\prod_{j=1}^l \sigma(t_j)}\exp \left(-x\sum_{i=1}^l \zeta(t_j)\right),\nonumber 
\end{align}
where $\sigma(x)$ is the Weierstrass sigma function and $\sigma_i(x)$ $(i=1,2,3)$ are the Weierstrass co-sigma functions. 
Consider fixing the value $E$. If $\Lambda (x,E) \not \in {\mathcal F}$, then $\Lambda (x,E)$ and $\Lambda (-x,E)$ are linearly independent and the linear combinations of $\Lambda (x,E)$ and $\Lambda (-x,E)$ exhaust solutions to equation (\ref{InoEF}).
Hence if there is no solution to (\ref{InoEF})  in ${\mathcal F}$, a basis of the solution is written as $\Lambda (x,E)$ and $\Lambda (-x,E)$ where $\Lambda (x,E)$ is given as (\ref{BA}). Note that the condition for $E$ that there is no solution to (\ref{InoEF})  in ${\mathcal F}$ is equivalent to the condition that the value $E$ satisfies $Q(E)\neq 0$.

Conversely assume that differential equation (\ref{InoEF}) has the solution written as 
\begin{equation}
\tilde{\Lambda }(x)=\frac{\prod_{j=1}^l \sigma(x+t_j)}{\sigma(x)^{l_0}\sigma_1(x)^{l_1}\sigma_2(x)^{l_2}\sigma_3(x)^{l_3}}\exp(cx),
\label{BV}
\end{equation}
for some $c$ and $t_j$ $(j=1,\dots ,l)$ such that $t_j\neq t_{j'}$ $(j\neq j')$ and $t_j \neq \omega _i$ $(j=1,\dots ,l, \: i=0,1,2,3)$.
By a straightforward calculation, it follows that
\begin{align}
& \frac{d^2\tilde{\Lambda }(x)}{dx^2}= \tilde{\Lambda }(x) \left( \sum_{i=0}^3 l_i(l_i +1)\wp(x+\omega_i)+ 2l_0 \zeta(x)\left(-c-\sum_{j=1}^l \zeta(t_j)\right) \right. \label{eqn:diff} \\
& + 2\sum_{i=1}^3 l_i\zeta(x+\omega_i)\left(-c-l\zeta(\omega_i)-\sum_{j=1}^l\zeta(t_j-\omega_i)\right) \nonumber \\
& + 2\sum_{j=1}^l \zeta(x+t_j)\left(c+\sum_{k\neq j} \zeta(t_k-t_j)-l_0\zeta(-t_j)-\sum_{i=1}^3l_i(\zeta(-t_j+\omega _i )-\zeta(\omega_i))\right)\nonumber \\
& +\sum_{j=1}^l\sum_{i=0}^3 l_i\left(\wp(t_j-\omega_i)-\zeta(t_j-\omega_i)^2\right)-\sum_{j<k}\left(\wp(t_j-t_k)-\zeta(t_j-t_k)^2\right) \nonumber \\
& \left. + c^2-(l_0l_1+l_2l_3)e_1-(l_0l_2+l_1l_3)e_2-(l_0l_3+l_1l_2)e_3+\sum_{i=1}^3l_i\eta_i (2c+l\eta_i) \right) . \nonumber 
\end{align}
Note that we used the formula
\begin{equation}
(\zeta(x+a)-\zeta(x+b)-\zeta(a-b))^2=\wp (x+a)+\wp(x+b)+\wp (a-b)
\end{equation}
to derive equation (\ref{eqn:diff}).

Hence we find that the function $\tilde{\Lambda }(x)$ satisfies equation (\ref{InoEF}) if and only if  $t_j$ $(j=1,\dots ,l)$ and $c$ satisfy the following equations;
\begin{align}
& \sum_{k\neq j} \zeta(-t_j+t_k)-l_0\zeta(-t_j)-\sum_{i=1}^3l_i(\zeta(-t_j+\omega _i )-\zeta(\omega_i))= -c\; \; (j=1, \dots ,l),
\label{BAeq}\\
& (1-\delta_{l_0,0})(c+\sum_{j=1}^l \zeta(t_j) )=0, \nonumber \\
& (1-\delta_{l_i,0})(c+l\zeta(\omega_i)+\sum_{j=1}^l \zeta(-\omega_i+t_j))=0, \; \; (i=1,2,3). \nonumber
\end{align}
The eigenvalue $E$ is given by
\begin{align}
& E=-c^2+(l_0l_1+l_2l_3)e_1+(l_0l_2+l_1l_3)e_2+(l_0l_3+l_1l_2)e_3-\sum_{i=1}^3l_i\eta_i (2c+l\eta_i) \label{BAE} \\
& -\sum_{j=1}^l\sum_{i=0}^3 l_i(\wp(t_j-\omega_i)-\zeta(t_j-\omega_i)^2)+\sum_{j<k}(\wp(t_j-t_k)-\zeta(t_j-t_k)^2). \nonumber 
\end{align}
Note that there are other expressions of $E$ in terms of $t_j$ $(j=1,\dots ,l)$ and $c$.

Summarizing, we obtain the following theorem:
\begin{thm} \label{BAthm}
Set 
\begin{equation}
\tilde{\Lambda }(x)=\frac{\prod_{j=1}^l \sigma(x+t_j)}{\sigma(x)^{l_0}\sigma_1(x)^{l_1}\sigma_2(x)^{l_2}\sigma_3(x)^{l_3}}\exp(cx), \; \; (l=l_0+l_1+l_2+l_3)
\label{BV2}
\end{equation}
and assume $t_j\neq t_{j'}$ $(j\neq j')$ and $t_j \neq \omega _i$ $(j=1,\dots ,l,\; i=0,1,2,3)$.

The condition that the function $\tilde{\Lambda }(x)$ satisfies equation (\ref{InoEF}) for some $E$ is equivalent to $t_j$ $(j=1,\dots ,l)$ and $c$ satisfying relations (\ref{BAeq}), with the eigenvalue $E$ given as (\ref{BAE}).
For the eigenvalue $E$ such that $Q(E)\neq 0$, a basis of solutions to (\ref{InoEF}) is written as $\tilde{\Lambda }(x)$ and $\tilde{\Lambda }(-x)$, and $t_j$ $(j=1,\dots ,l)$ such that $t_j\neq t_{j'}$ $(j\neq j')$ and $t_j \neq \omega _i$ $(j=1,\dots ,l,\; i=0,1,2,3)$.
\end{thm}
\begin{remk}
(i) Comparing with \cite{FVthr,Tak}, it would be reasonable to call the function $\tilde{\Lambda }(x)$ the Bethe vector and equation (\ref{BAeq}) the Bethe Ansatz equation.

(ii) For each $l_0,l_1,l_2,l_3$, the number of eigenvalues $E$ which do not have a basis of solutions to (\ref{InoEF}) as the Bethe vectors (\ref{BV2}) is at most the dimension of the space $\tilde{V}$ (see (\ref{Vdecom})) which is given in Theorem \ref{periodic}.

(iii) If $l_0 \neq 0$ and $l_1 \neq 0$ then the number of equations (\ref{BAeq}) is strictly greater than the number of the variables $t_j$ $(j=1,\dots ,l)$, although equations (\ref{BAeq}) have solutions $(t_1, \dots ,t_l)$ if $Q(E)\neq 0$.

For the case $l_1=l_2=l_3=0$ (the case of the Lam\'e equation), the equation $c+\sum_{j=1}^l \zeta(t_j) =0$ is obtained by summing up other equations in (\ref{BAeq}). Hence the number of equations is essentially equal to that of the variables $t_j$ $(j=1,\dots ,l_0)$.
\end{remk}

\subsection{Bethe Ansatz equation in terms of theta functions} \label{sec:BAtheta}

We will introduce the Bethe Ansatz equation written in terms of theta functions.

Set $\tau = \omega_3/\omega_1$ and  $\omega_1 =1/2$, then the Hamiltonian (see (\ref{Ino})) is described as 
\begin{align}
& H= -\frac{d^2}{dx^2} +l_0(l_0+1)\wp (x)\label{Hamil} \\
& +l_1(l_1+1)\wp \left(x+\frac{1}{2}\right)+l_2(l_2+1)\wp \left(x+\frac{1+\tau}{2}\right)+l_3(l_3+1)\wp \left(x+\frac{\tau}{2}\right). \nonumber
\end{align}

By replacing the sigma function with the theta function by formula (\ref{sigthetaA}), it is seen that if $Q(E)\neq 0$, then the basis of solutions to (\ref{InoEF}) is written as $\tilde{\Lambda }(x)$ and $\tilde{\Lambda }(-x)$, where
\begin{equation}
\tilde{\Lambda }(x)=\frac{\exp(\pi \sqrt{-1}cx)\prod_{j=1}^l \theta(x+t_j)}{\theta(x)^{l_0}\theta(x+1/2)^{l_1}\theta(x+(1+\tau)/2)^{l_2}\theta(x+\tau/2)^{l_3}},
\label{tBV}
\end{equation}
$t_j \neq t_{j'}$ $(j\neq j')$, $t_j \not \in \frac{1}{2} \Zint + \frac{\tau}{2} \Zint $ $(j=1, \dots ,l)$ and $l=l_0+l_1+l_2+l_3$.

By similar calculations as the previous subsection, we obtain relations for $t_j$ $(j=1,\dots ,l)$ and $c$ for which the function $\tilde{\Lambda }(x)$ becomes an eigenfunction of operator (\ref{Hamil}).

\begin{thm} \label{thm:thetaBA}
The function $\tilde{\Lambda }(x)$ (see (\ref{tBV})) with the condition $t_j\neq t_{j'}$ $(j\neq j')$ and $t_j \not \in \frac{1}{2} \Zint + \frac{\tau}{2} \Zint $ is an eigenfunction of (\ref{Hamil}) if and only if  $t_j$ $(j=1,\dots ,l(=l_0+l_1+l_2+l_3))$ and $c$ satisfy the relations,
\begin{align}
& c\pi\sqrt{-1} +\sum_{k\neq j}\frac{\theta '(t_k-t_j)}{\theta (t_k-t_j)}+l_0\frac{\theta '(t_j)}{\theta (t_j)}+l_1\frac{\theta '(t_j-1/2)}{\theta (t_j-1/2)} \label{thetabae}\\
& +l_2\frac{\theta '(t_j-(1+\tau)/2)}{\theta (t_j-(1+\tau)/2)}+l_3\frac{\theta '(t_j-\tau/2)}{\theta (t_j-\tau/2)}=0, \; \; \; (j=1,\dots ,l) \nonumber \\
& (1-\delta_{l_0,0})\left(c\pi\sqrt{-1} + \pi \sqrt{-1}(l_2+l_3)+ \sum_{j=1}^l\frac{\theta '(t_j)}{\theta (t_j)}\right) =0, \nonumber \\
& (1-\delta_{l_1,0})\left( c\pi\sqrt{-1} +\pi \sqrt{-1}(l_2+l_3)+ \sum_{j=1}^l\frac{\theta '(t_j-1/2)}{\theta (t_j-1/2)}\right) =0, \nonumber \\
& (1-\delta_{l_2,0})\left(c\pi\sqrt{-1}  -\pi \sqrt{-1}(l_0+l_1)+ \sum_{j=1}^l\frac{\theta '(t_j-(1+\tau)/2)}{\theta (t_j-(1+\tau)/2)}\right) =0, \nonumber \\
& (1-\delta_{l_3,0})\left(c\pi\sqrt{-1}  -\pi \sqrt{-1}(l_0+l_1)+ \sum_{j=1}^l\frac{\theta '(t_j-\tau/2)}{\theta (t_j-\tau/2)}\right) =0. \nonumber 
\end{align}
The eigenvalue $E$ is given by
\begin{align}
& E=\pi^2(c^2+(l_0+l_1)(l_2+l_3))+l(l+1)\eta_1 
-\sum_{j<k}\frac{\theta ''(t_j-t_k)}{\theta (t_j-t_k)} \label{thetaE} \\
& +(l_0l_1+l_2l_3)e_1+(l_0l_2+l_1l_3)e_2+(l_0l_3+l_1l_2)e_3 \nonumber \\
& +\sum_{j=1}^l\left( l_0 \frac{\theta ''(t_j)}{\theta (t_j)}+l_1\frac{\theta ''(t_j-1/2)}{\theta (t_j-1/2)}+l_2\frac{\theta ''(t_j-(1+\tau)/2)}{\theta (t_j-(1+\tau)/2)}+l_3\frac{\theta ''(t_j-\tau/2)}{\theta (t_j-\tau/2)}\right) . \nonumber 
\end{align}

If the value $E$ satisfies $Q(E)\neq 0$,
then a basis of solutions to (\ref{InoEF}) is written as $\tilde{\Lambda}(x)$ and $\tilde{\Lambda}(-x)$, where the function $\tilde{\Lambda}(x)$ is given in (\ref{tBV}) and satisfies the condition $t_j\neq t_{j'}$ $(j\neq j')$ and $t_j \not \in \frac{1}{2} \Zint + \frac{\tau}{2} \Zint $ $(j=1, \dots ,l)$. 
\end{thm}

It would be hopeless to attempt to solve equations (\ref{thetabae}) explicitly for generic $\tau $.
Instead, it may be possible to solve them for the case $\tau \rightarrow \sqrt{-1}\infty$ and to gain some knowledge for the case $\tau \neq \sqrt{-1}\infty$ from this limit.
We will investigate solutions to equations (\ref{thetabae}) for the trigonometric $(\tau \rightarrow\sqrt{-1} \infty)$ case in the next section.

\section{Trigonometric limit} \label{sec:trig}
\subsection{Trigonometric limit of the polynomial $Q(E)$} \label{sec:trigPE}

In this section, we will consider the trigonometric limit $\tau \rightarrow \sqrt{-1}\infty$.

The trigonometric limits of the elliptic functions are given by
\begin{align}
& \wp(x) \rightarrow -\frac{\pi^2}{3}+\frac{\pi^2}{(\sin \pi x)^2}, \; \; 
\wp\left(x+\frac{1}{2} \right)  \rightarrow -\frac{\pi^2}{3}+\frac{\pi^2}{(\cos \pi x)^2}, \\
&  \wp\left(x+\frac{\tau}{2}\right)  \rightarrow -\frac{\pi^2}{3}, \; \; \wp\left(x+\frac{1+\tau}{2}\right)  \rightarrow -\frac{\pi^2}{3}, \; \; \frac{\theta(x)}{\theta'(0)} \rightarrow \sin \pi x,
\end{align}
as $\tau \rightarrow \sqrt{-1}\infty$. They converge uniformly on compact sets.
Due to this limit, the term $l_2(l_2+1)\wp(x+\frac{1+\tau}{2})+ l_3(l_3+1)\wp(x+\frac{\tau}{2}) $ converges to a constant. From now on, we consider the case $l_2=l_3=0$.
Set $p=\exp (2 \pi \sqrt{-1}\tau)$.
Note that the limit $\tau \rightarrow \sqrt{-1}\infty$ corresponds to the limit $p\rightarrow 0$. 

As $p\rightarrow 0$, we obtain $H \rightarrow  H_T-C_T$, where
\begin{align}
& H_T= -\frac{d^2}{dx^2} +l_0(l_0+1)\frac{\pi^2}{(\sin \pi x)^2}
+l_1(l_1+1)\frac{\pi^2}{(\cos \pi x)^2},
\label{trigH} \\
& C_T=(l_0(l_0+1)+l_1(l_1+1))\pi^2/3. \label{trigCT}
\end{align}

We will consider the trigonometric limit of functions which have appeared in section \ref{sec:monod}. Note that $e_1\rightarrow \frac{2\pi^2}{3}$, $e_2\rightarrow \frac{-\pi^2}{3}$ and $e_3\rightarrow \frac{-\pi^2}{3}$ as $p\rightarrow 0$.
The function $\Xi (x,E)$ (see (\ref{Fx})) admits the trigonometric limit 
\begin{equation} 
\Xi _T(x,E)=\tilde{c}_0(E)+ \sum_{j=0}^{l_0-1} \frac{\tilde{a}^{(0)}_j(E)}{(\sin\pi x)^{2(l_0-j)}}+\sum_{j=0}^{l_1-1} \frac{\tilde{a}^{(1)}_j(E)}{(\cos \pi x)^{2(l_1-j)}},
\label{FTx}
\end{equation}
and the function $\Xi _T(x,E)$ satisfies the trigonometric version of equations (\ref{prodDE}) and (\ref{const}). We will calculate the polynomial $Q_T(E)$ which is the trigonometric version of $Q(E)$. From the trigonometric version of equation (\ref{const}), we have $Q_T(E)=(E+C_T) \tilde{c}_0(E)^2$.

\begin{prop} \label{con:triconst}
If $l_0+ l_1$ is even, then
\begin{equation}
Q_T(E)=
b\cdot (E+C_T)\prod_{i=1}^{\frac{l_0+l_1}{2}}(E-(2i)^2+C_T) ^2\prod_{i=1}^{\frac{|l_0-l_1|}{2}}(E-(2i-1)^2+C_T)^2,
\end{equation}
where $C_T$ is defined in (\ref{trigCT}) and $b$ is a constant which is independent of $E$.

If $l_0+ l_1$ is odd, then
\begin{equation}
Q_T(E)=b\cdot(E+C_T)\prod_{i=1}^{\frac{|l_0-l_1|-1}{2}}(E-(2i)^2+C_T)^2 \prod_{i=1}^{\frac{l_0+l_1+1}{2}}(E-(2i-1)^2+C_T)^2.
\end{equation}
\end{prop}
\begin{proof}
For simplicity, we prove for the case $l_0\geq l_1$ and $l_0, l_1 \in 2\Zint _{>0}$. The proofs for the other cases are similar.

From the proof of Proposition \ref{propconj1}, it follows that the polynomial $Q_T(E)$ is obtained by a suitable limit of the characteristic polynomial of the space $\tilde{V}$ (see (\ref{Vdecom})). For this case, we have
$\tilde{V}= \tilde{V}_{1,1,1,1} \oplus  \tilde{V}_{1,1,-1,-1} \oplus  \tilde{V}_{1,-1,-1,1} \oplus  \tilde{V}_{1,-1,1,-1}$.

Set $\epsilon =(e_2-e_3)$, then $\epsilon \rightarrow 0$ as $\tau \rightarrow \sqrt{-1} \infty$. We will solve recursive relation (\ref{rec:invsp}) for the case $\epsilon \rightarrow 0$. 
If $\epsilon \neq 0$ then the invariant subspace $\tilde{V}_{1,1,1,1}$ is given as $\tilde{V}_{1,1,1,1}= \bigoplus _{n=0}^{\frac{l_0+l_1}{2}} \Cplx \left( \frac{\sigma_1(x)}{\sigma (x)} \right) ^{-\frac{l_1}{2}} \wp(x)^n$. Then the characteristic polynomial of the space  $\tilde{V}_{1,1,1,1}$ for $\epsilon \rightarrow 0$ is expressed as
\begin{equation}
\epsilon ^{-\frac{l_0+l_1}{2}-1} c_1\prod_{i=0}^{\frac{l_0+l_1}{2}}(E-(2i)^2+C_T) + \epsilon ^{-\frac{l_0+l_1}{2}}\star +\dots,
\end{equation}
where $c_1$ is a constant. Similarly, the leading terms of the characteristic polynomials of the spaces $\tilde{V}_{1,1,-1,-1}$, $\tilde{V}_{1,-1,-1,1}$, $\tilde{V}_{1,-1,1,-1}$ as $\epsilon \rightarrow 0$ are given as
\begin{equation} 
\prod_{i=1}^{\frac{l_0+l_1}{2}}\left(E-(2i)^2+C_T\right)  ,\; 
\prod_{i=1}^{\frac{l_0-l_1}{2}}\left(E-(2i-1)^2+C_T\right) , \; 
\prod_{i=1}^{\frac{l_0-l_1}{2}}\left(E-(2i-1)^2+C_T\right) ,
\end{equation}
up to constant multiplications respectively.

By multiplying the leading terms of these four cases, we obtain the proposition.
\end{proof}

\subsection{Bethe Ansatz equation for the trigonometric model}

In this subsection, we will obtain the trigonometric version of equations (\ref{thetabae}).

Set $X=\exp(2\pi \sqrt{-1}x)$ and $T_j=\exp(2\pi \sqrt{-1} t_j)$. 
Now we state the trigonometric version of Theorem \ref{thm:thetaBA}.
Set 
\begin{equation}
\tilde{\Lambda }_T(X)=\frac{\prod_{j=1}^{l_0+l_1} (1-XT_j)}{(X-1)^{l_0}(X+1)^{l_1}}X^{c/2},
\label{triBV}
\end{equation}
and assume $T_j\neq T_{j'}$ $(j\neq j')$ and $T_j \neq 0, \pm 1$ $(j\in \{1,\dots ,l_0+l_1 \})$.
Suppose that $T_j$ and $c$ satisfy relations,
\begin{align}
& c+\sum_{k\neq j}\frac{T_k+T_j}{T_k-T_j}+l_0\frac{T_j+1}{T_j-1}+l_1\frac{T_j-1}{T_j+1}=0 \; \; \; (j=1,\dots ,l_0+l_1) \label{tBAE1} \\
& (1-\delta_{l_0,0})\left(c+\sum_{j=1}^{l_0+l_1}\frac{T_j+1}{T_j-1} \right)=0 , \; \; \; (1-\delta_{l_1,0})\left(c+\sum_{j=1}^{l_0+l_1}\frac{T_j-1}{T_j+1} \right)=0 \label{tBAE2},
\end{align}
then it is shown that the function $\tilde{\Lambda }_T(e^{2\pi\sqrt{-1}x})$ is an eigenfunction of operator $H_T$ (\ref{trigH}) with an eigenvalue $E=\pi ^2 c^2$.

Conversely, if the function $\tilde{\Lambda }_T(X)=\frac{\prod_{j=1}^{l_0+l_1} (1-XT_j)}{(X-1)^{l_0}(X+1)^{l_1}}X^{c/2}$ ($T_j\neq T_{j'}$ ($j\neq j'$) and $T_j\neq 0, \pm 1$ ($j=1,\dots ,l_0+l_1$)) satisfies the equation $H_T \tilde{\Lambda }_T(e^{2\pi\sqrt{-1}x})= E\tilde{\Lambda }_T(e^{2\pi\sqrt{-1}x})$, we obtain equations (\ref{tBAE1}, \ref{tBAE2}) and $E=\pi ^2 c^2$.

\begin{prop} \label{prop:comp}
Fix the value $E$. If $Q_T(E)\neq 0$ then the basis of the solutions to the equation $H_T f(x)=(E+C_T)f(x)$ is written as $\tilde{\Lambda }_T(e^{2\pi\sqrt{-1}x})$ and $\tilde{\Lambda }_T(e^{-2\pi\sqrt{-1}x})$ where $\tilde{\Lambda }_T(e^{2\pi\sqrt{-1}x})$ is given in (\ref{triBV}) with the conditions $c^2=\frac{E+C_T}{\pi^2}$, $T_j\neq T_{j'}$ $(j\neq j')$ and $T_j \neq 0, \pm 1$ $(j\in \{1,\dots ,l_0+l_1 \})$.
\end{prop}
\begin{proof}
From expression (\ref{FTx}), it follows that
\begin{align} 
& \Xi _T(x,E)=\frac{\tilde{c}_0(E)(\sin\pi x)^{2(l_0+l_1)}+ (\mbox{lower terms in } (\sin\pi x)^2)}{(\sin\pi x)^{2l_0}(\cos \pi x)^{2l_1}} ,\label{FTx2} 
\end{align}
for some $t_1,\dots ,t_{l_0+l_1}$, where the polynomial $\tilde{c}_0(E)$ is common in (\ref{FTx}) and (\ref{FTx2}). From the relation $Q_T(E)=(E+C_T) \tilde{c}_0(E)^2$, we have $\tilde{c}_0(E) \neq 0$. Hence $\Xi _T(x,E)$ is written as
\begin{align} 
& \Xi _T(x,E)=\frac{\tilde{c}_0(E)\prod_{k=1}^{l_0+l_1}((\sin\pi x)^{2}- (\sin\pi t_j)^{2})}{(\sin\pi x)^{2l_0}(\cos \pi x)^{2l_1}}
\label{FTx3}
\end{align}
for some $t_1,\dots ,t_{l_0+l_1}$. 
By a similar argument to that performed in section \ref{sec:BAEsigma}, we obtain
\begin{equation}
\Lambda _T(x)= \frac{\sqrt{\tilde{c}_0(E)}\prod_{k=1}^{l_0+l_1}\sin \pi (x+t_k)}{(\sin\pi x)^{l_0}(\cos \pi x)^{l_1}}\exp(\pi \sqrt{-1}c),
\end{equation}
where the value $c$ satisfies $c^2=\frac{E+C_T}{\pi^2}$ and the signs of $t_1,\dots ,t_{l_0+l_1}$ are suitably chosen.

If $\Lambda _T(x)$ satisfies the differential equation $(H_T -C_T -E)\Lambda _T(x)=0$, $\Lambda _T(-x)$ also satisfies it. From the trigonometric version of formula (\ref{integ}), it is seen that the functions $\Lambda _T(x)$ and $\Lambda _T(-x)$ are linearly independent. By a similar argument to the proof of Proposition \ref{prop:tjdist}, the values $ \pm t_j$, $(j=1,\dots ,l_0+l_1)$, $0, \; 1/2$ are mutually different.
By setting $X=\exp(2\pi \sqrt{-1}x)$ and $T_j=\exp (2\pi \sqrt{-1} t_j)$, we obtain the proposition.
\end{proof}

\subsection{Solutions to the Bethe Ansatz equation}

We will solve the trigonometric Bethe Ansatz equation.

\begin{thm}
Assume $Q_T(\pi ^2c^2-C_T)\neq 0$.
Let $(T_{1}^0, \dots ,T_{l}^0)$ $(l=l_0+l_1)$ be a solution to equations (\ref{tBAE1}, \ref{tBAE2}) and let $\sigma_i=\sum_{j_1<\dots<j_i}T_{j_1}^0\dots T_{j_i}^0$ be the $i$-th elementary symmetric function of the solution $(T_{1}^0, \dots ,T_{l}^0)$.
Then there exists a solution to (\ref{tBAE1}, \ref{tBAE2}) uniquely up to permutation of  $(T_{1}^0, \dots ,T_{l}^0)$. The functions $\sigma_i$ are determined by the following recursive relation,
\begin{align}
& \sigma_0=1, \; \; \; \; \sigma_1=\frac{(l_0-l_1)(c-l_0-l_1)}{c+1}, \\
&  \sigma_i =\frac{(l_0-l_1)(c+2i-2-l_0-l_1)}{i(i+c)}\sigma_{i-1}+ \frac{(i-2-l_0-l_1)(c+i-2-l_0-l_1)}{i(i+c)} \sigma_{i-2} \label{sigmarecrel}
\end{align}
for $i=2,3,\dots ,l_0+l_1$.
\end{thm}
\begin{proof}
By the assumption $Q_T(\pi ^2c^2-C_T)\neq 0$ and Proposition \ref{prop:comp}, there exists a solution to the equation $(H_T-\pi^2c^2)\tilde{\Lambda }_T(e^{2\pi \sqrt{-1} x})=0$ of the form $\tilde{\Lambda }_T(X)=\frac{\prod_{j=1}^l (1-XT_j)}{(X-1)^{l_0}(X+1)^{l_1}}X^{c/2}$. Then the values $(T_1,\dots, T_l)$ satisfy equations (\ref{tBAE1}, \ref{tBAE2}) and we have
\begin{equation}
\tilde{\Lambda }_T(X)=\frac{1+\sum_{j=1}^l(-1)^j\sigma_j X^j}{(X-1)^{l_0}(X+1)^{l_1}}X^{c/2}.
\end{equation}
By comparing the coefficients of $X^j$ $(j=1,\dots ,l)$ of the equation $(H_T-\pi^2c^2)\tilde{\Lambda }_T(X)=0$ ($X=e^{2\pi \sqrt{-1} x}$), we obtain recursive relation (\ref{sigmarecrel}) for $\sigma_i$.
\end{proof}

Now we exhibit explicit expressions of $\sigma_m$.
By straightforward calculations, we obtain the following proposition:
\begin{prop} \label{prop:exp}

(i) If $l_1=0 $ then
\begin{equation}
\sigma_{m}=\prod_{j=0}^{m-1}\frac{(l_0-j)(c+j-l_0)}{(j+1)(c+j+1)}.
\end{equation}

(ii) If $l_1=1 $ then
\begin{align}
& \sigma_{m}=\frac{\prod_{j=0}^{m-2}(l_0-j)\prod_{j=1}^{m-2}(c+j-l_0)}{m!\prod_{j=0}^{m-1}(c+j+1)} \cdot \\
& \; \; \; \; \; \cdot (c-l_0-1)((l_0-2m+1)c-l_0^2-l_0+2ml_0+2m-2m^2) \nonumber
\end{align}

(iii) If $l_1=2 $ then
\begin{align}
& \sigma_{m}=\frac{\prod_{j=0}^{m-3}(l_0-j)\prod_{j=2}^{m-3}(c+j-l_0)}{m!\prod_{j=0}^{m-1}(c+j+1)}(c-l_0)(c-l_0-2)\cdot \\
& \; \; \; \cdot (c^2(l_0^2-(4m-3)l_0+4m^2-8m+2)\nonumber \\
& \; \; \; \; \; \; +2c(-l_0^3+(4m-3)l_0^2  -(6m^2-10m+2)l_0+4m^3-12m^2+8m) \nonumber \\
& \; \; \; \; \; \;+l_0^4-(4m-3)l_0^3 +(8m^2-12m+1)l_0^2 \nonumber \\
& \; \; \;  \; \; \; -(8m^3-20m^2+8m+3)l_0+4m^4-16m^3+16m^2-2). \nonumber
\end{align}

(iv) If $l_1=l_0$ then
\begin{equation}
 \sigma_{2m-1}=0 , \; \; \; \;  \sigma_{2m}=\prod_{j=0}^{m-1}\frac{(j-l_0)(c+2j-2l_0)}{(j+1)(c+2j+2)}.
\end{equation}

(v) If $l_1=l_0 -1$ then
\begin{align}
& \sigma_{2m-1}=\frac{\prod_{j=1}^{m-1}(j-l_0)\prod_{j=0}^{m-1}(c+2j+1-2l_0)}{(m-1)!\prod_{j=0}^{m-1}(c+2j+1)}, \\
& \sigma_{2m}=\frac{\prod_{j=1}^{m}(j-l_0)\prod_{j=0}^{m-1}(c+2j+1-2l_0)}{m!\prod_{j=0}^{m-1}(c+2j+1)}.\nonumber
\end{align}

(vi)  If $l_1=l_0 -2$ then
\begin{align}
& \sigma_{2m-1}=\frac{2\prod_{j=2}^{m}(j-l_0)\prod_{j=0}^{m-1}(c+2j+2-2l_0)}{(m-1)!(c+1)\prod_{j=0}^{m-2}(c+2j+2)}, \\
& \sigma_{2m}=-\frac{\prod_{j=2}^{m}(j-l_0)\prod_{j=0}^{m-1}(c+2j+2-2l_0)}{m!(c+1)\prod_{j=0}^{m-1}(c+2j+2)}\cdot \nonumber \\
& \; \; \; \; \; \cdot (c(l_0-2m-1)+(4m+1)l_0-(2m+1)^2).\nonumber
\end{align}
\end{prop}

\begin{remk}
Proposition \ref{prop:exp} (i) reproduces Proposition 4.1 in \cite{Tak}.
\end{remk}

\section{Eigenstates for the $BC_1$ Inozemchev system} \label{sec:pert}
\subsection{Eigenstates for the trigonometric Hamiltonian}

Let us find square-integrable eigenstates of the Hamiltonian $H_T$.
Set
\begin{equation}
\Phi(x)=(\sin \pi x)^{l_0+1}(\cos \pi x)^{l_1+1}, \; \; \; {\mathcal H}_T=\Phi(x)^{-1} H_T \Phi(x),
\end{equation}
then the gauge transformed Hamiltonian is written as
\begin{equation}
{\mathcal H}_T= -\frac{d^2}{dx^2}-2\pi \left(
\frac{(l_0 +1)\cos \pi x}{\sin \pi x}- \frac{(l_1 +1)\sin \pi x}{\cos \pi x}
\right) \frac{d}{dx}+(l_0+l_1+2)^2\pi^2.
\end{equation}
By a straightforward calculation we obtain
\begin{align}
& {\mathcal H}_T (\cos2\pi x)^m=\pi^2(2m+l_0+l_1+2)^2(\cos2\pi x)^m \\
& \; \; \; \; \; \; +4\pi^2((-m(m-1)+m(l_0-l_1))(\cos2\pi x)^{m-1}. \nonumber
\end{align}
Therefore the operator ${\mathcal H}_T$ acts triangularly on the space spanned by $\{ 1, \cos2\pi x, (\cos 2\pi x)^2 ,\dots \}$. Since $(2m+l_0+l_1+2)^2 \neq (2m'+l_0+l_1+2)^2$ $(m, m' \in \Zint _{\geq 0}, \; m\neq m')$, the operator ${\mathcal H}_T$ is diagonalizable and
\begin{equation}
{\mathcal H}_T \psi_m(x)=\pi^2(2m+l_0+l_1+2)^2\psi_m(x),
\end{equation}
for some functions $\psi_m(x)=(\cos2\pi x)^m+ $(lower terms). Set $w=(1-\cos 2\pi x)/2$. By a change of variable, it follows that
\begin{align}
& {\mathcal H}_T -\pi^2(2m+l_0+l_1+2)^2= -4\pi ^2 \cdot \label{Ghyper}  \\
& \left\{ w(1-w)\frac{d^2}{dw^2}+\left(\frac{2l_0+3}{2}-((l_0+l_1+2)+1)w\right)\frac{d}{dw}+m(m+l_0+l_1+2) \right\}. \nonumber
\end{align}
Hence the operator ${\mathcal H}_T -\pi^2(2m+l_0+l_1+2)^2$ is transformed to the hypergeometric operator. From expression (\ref{Ghyper}), we obtain
\begin{equation}
\psi_m(x) = \tilde{c}_m G_m\left( l_0+l_1+2, \frac{2l_0+3}{2}; \frac{1-\cos 2\pi x}{2}\right),
\label{Jacobipol}
\end{equation}
where the function $ G_m(\alpha, \beta; w)=_2 \! F_1(-m,\alpha+m;\beta;w)$ is the Jacobi polynomial of degree $m$ and $\tilde{c}_m$ is a constant.

We define the inner product by 
\begin{equation}
\langle f,g\rangle =\int_{0}^{1} dx\overline{f(x)} g(x),
\end{equation}
then the values $ \langle \Phi(x)\psi_m(x),\Phi(x)\psi_{m'}(x)\rangle$ are finite and
\begin{align}
& \pi^2(2m+l_0+l_1+2)^2\langle \Phi(x)\psi_m(x),\Phi(x)\psi_{m'}(x)\rangle\\
& =\langle H_T \Phi(x)\psi_m(x),\Phi(x)\psi_{m'}(x)\rangle 
= \langle \Phi(x)\psi_m(x),H_T\Phi(x)\psi_{m'}(x)\rangle\nonumber \\
& =\pi^2(2m'+l_0+l_1+2)^2\langle \Phi(x)\psi_m(x),\Phi(x)\psi_{m'}(x)\rangle . \nonumber
\end{align}
Therefore the functions $\Phi(x)\psi_m(x)$ $(m\in \Zint_{\geq 0})$ are orthogonal.

\begin{remk}
(i) The values $\tilde{c}_m$ and  $ \langle \Phi(x)\psi_m(x),\Phi(x)\psi_{m}(x)\rangle$ are known explicitly. \\
(ii) The set of functions $\Phi(x)\psi_m(x)$ $(m\in \Zint_{\geq 0})$ is a complete orthogonal basis of the space $L^2(\Rea /\Zint )$. For the proof, see \cite{KT} or \cite{Tak2}.\\
(iii) The functions $\psi_m(x)$ are simply the $BC_1$-Jacobi polynomials.
\end{remk}

\subsection{Jacobi polynomial, Bethe Ansatz and perturbation}

In this subsection, we will see the relationship between the Jacobi polynomial and the function obtained by the Bethe Ansatz method. After that, we will show the holomorphy of the perturbation with respect to the parameter $p=\exp (2\pi \sqrt{-1} \tau) \sim 0$ for the Hamiltonian of the Inozemtsev model (\ref{Ino}) with the condition $l_2=l_3=0$. For the case $l_1=l_2=l_3=0$, the corresponding result is obtained in \cite{Tak}.

Let us recall the Bethe Ansatz method for the trigonometric model, which was discussed in section \ref{sec:trig}.

Let $m$ be an non-negative integer, and set $c=l_0+l_1+2+2m$, then $Q_T(\pi^2 c^2 - C_T)\neq 0$, where $Q_T(E)$ and $C_T$ were defined in section \ref{sec:trigPE}. Hence the solution $T_1,\dots ,T_{l_0+l_1}$ of the Bethe Ansatz equation (see (\ref{tBAE1}, \ref{tBAE2})) with the condition $c=l_0+l_1+2+2m$ exists uniquely up to permutation. In this setting, the function 
\begin{equation}
\Lambda _T(x)=\frac{\prod_{j=1}^{l_0+l_1}\sin \pi (x+t_j)}{(\sin \pi x)^{l_0} (\cos \pi x)^{l_1}}e^{\pi \sqrt{-1} cx},
\end{equation}
satisfies $H_T\Lambda _T(x)= \pi^2 c^2 \Lambda _T(x)$, where $\exp(2\pi \sqrt{-1} t_j)=T_j$ and the values $(T_1,\dots ,T_{l_0+l_1})$ are a solution to equations (\ref{tBAE1}, \ref{tBAE2}). From the condition $Q_T(\pi^2 (l_0+l_1+2+2m)^2 - C_T)\neq 0$, it follows that the functions $\Lambda _T(x)$ and $\Lambda _T(-x)$ are linearly independent.
Set $\Lambda^{sym}_T(x)=\Lambda _T(x) -(-1)^{l_0}\Lambda _T(-x)$, then $\Lambda^{sym}_T (x) \neq 0$ and the degree of the pole of the function $\Lambda^{sym}_T(x)$ at $x=n$ (resp. $x=n+1/2$) $(n\in \Zint )$ is less than $l_0$ (resp. $l_1$). Since the exponents of the differential equation $(H_T -\pi^2 c^2 )\Lambda^{sym}_T(x) =0$ at $x=n$ (resp. $x=n+1/2$) are $-l_0$ and $l_0+1$ (resp. $-l_1$ and $l_1+1$), the function $\Lambda^{sym}_T(x)$ has zeros of degree $l_0+1$ (resp. $l_1+1$) at $x=n$ (resp. $x=n+1/2$).
Since the eigenvalue of the function $\Phi(x) \psi _m(x)$ with respect to the Hamiltonian $H_T$ is $\pi ^2 (l_0+l_1+2+2m)^2$ and the function $\Phi(x) \psi _m(x)$ is holomorphic at $x=0$, it follows that $\Lambda^{sym}_T(x) =B_m \Phi(x)\psi _m(x)$ for some non-zero constant $B_m$.

Now, we will show that for each $c= l_0+l_1+2+2m$ $(m \in \Zint_{\geq 0})$ there exists $\tilde{p} \in \Rea _{>0}$ such that the solution to Bethe Ansatz equation (\ref{thetabae}) exists for each $p$ $(|p|<\tilde{p})$ and it converges  analytically to the solution to the Bethe Ansatz equation for the case $p=0$. Set $c_m=l_0+l_1+2+2m$.
Since the coefficients of the polynomial $Q(E)$ defined in section \ref{sec:monod} are holomorphic in $p$, there exists $\epsilon_1, p_1 \in \Rea _{>0} $ such that $Q(E) \neq 0$ for all $E$ and $p$ such that $|E+ C_T-\pi^2 c_m^2 |<\epsilon_1$ and $ |p|<p_1$.
Let us recall that the function $\Xi (x,E)$ is expressed as
\begin{equation}
\Xi (x,E)=\frac{\prod_{j=1}^{l_0+l_1} (\wp (x) -\wp(t_j))}{\wp(x)^{l_0} \wp(x+1/2)^{l_1}}
\label{contF}
\end{equation}
Set $a_j=\wp(t_j)$ and  $z=\wp (x)$, then the values $a_j$ depend on the values $E$ and $p$.

The numerator of equation (\ref{contF}) is written as a polynomial in $z$ and the coefficients of $z^j$ $(j=1, \dots ,l_0+l_1)$ are holomorphic in $E$ and $p$. Since the zeros of the numerator of (\ref{contF}) in $z$ are simple, the values $a_j$ are holomorphic in $E$ and $p$ such that $|E+ C_T-\pi^2 c_m^2 |<\epsilon_1, \; |p|<p_1$.

Suppose $|E+ C_T-\pi^2 c_m^2 |<\epsilon_1$ and $|p|<p_1$, then $Q(E)\neq 0$ and there exists a solution to Bethe Ansatz equation  (\ref{thetabae}).
Denote the solution by $(t_1(E,p), \dots ,t_{l_0+l_1}(E,p))$. If we choose the index suitably, then $a_j=\wp(t_j(E,p))$.

From Proposition \ref{prop:tjdist} and the proof of Proposition \ref{prop:comp}, it follows that the values $t_j(E,p)$ are mutually distinct and are not contained in  $\left\{ \begin{array}{ll}
\Zint /2 & (p=0)\\
\Zint /2+ \Zint\tau /2 & (p\neq 0) \end{array} \right. $ for $E$ and $p$ such that $|E+ C_T-\pi^2 c_m^2 |<\epsilon_1, \; |p|<p_1$.
Since $t_j(E,p) \not \in \frac{1}{2} \Zint +\frac{1}{2} \tau \Zint$, we have $\wp '(t_j(E,p)) \neq 0$ and the values $t_j(E,p)$ are holomorphic in $E$ and $p$ because the values $a_j$ are holomorphic in $E$ and $p$.
From (\ref{thetaE}), the value $c$ satisfies the relation $\pi^2c^2=h(E,p)$ where
\begin{align}
& h(E,p)=E-(l_0+l_1)(l_0+l_1+1)\eta_1 -l_0l_1e_1 \\
& +\sum_{j<k}\frac{\theta ''(t_j(E,p)-t_k(E,p))}{\theta (t_j(E,p)-t_k(E,p))} -\sum_{j=1}^{l_0+l_1}\left( l_0 \frac{\theta ''(t_j(E,p))}{\theta (t_j(E,p))}+l_1\frac{\theta ''(t_j(E,p)-1/2)}{\theta (t_j(E,p)-1/2)}\right) . \nonumber 
\end{align}

We are going to investigate the Bethe Ansatz equation with the condition $c=c_m$ for each $p$ ($|p|$: sufficiently small).
For the case $p=0$, we have $h(E-C_T,0)=E$ and $\frac{\partial h}{\partial E}(\pi^2c_m^2-C_T,0)=1$.
By the implicit function theorem, there exists $\tilde{p} \in \Rea_{>0}$ such that there exists a holomorphic function $E(p)$ in $p$ on the domain $|p|<\tilde{p}$ such that $h(E(p),p)=\pi^2 c_m^2$ and $E(0)=\pi^2 c_m^2-C_T$.
Set $\tilde{t}_j(p)=t_j(E(p),p)$ $(j=1,\dots ,l_0+l_1)$, then the function $\tilde{t}_j(p)$ is holomorphic in the variable $p$ on $|p|<\tilde{p}$, $(\tilde{t}_1(p), \dots , \tilde{t}_{l_0+l_1}(p))$ is a solution to Bethe Ansatz equation (\ref{thetabae}) with the condition $c=c_m$, and the corresponding eigenvalue $E(p)$ is also holomorphic in $p$.

Hence we obtain that the solution to Bethe Ansatz equation (\ref{thetabae}) at $p(=\exp(2\pi \sqrt{-1} \tau ))$ is holomorphically connected to the solution at $p= 0$ with the condition $c= c_m=l_0+l_1+2+2m$ $(m \in \Zint_{\geq 0})$.

Now, we set
\begin{align}
& \Lambda_m(x,p)= \frac{\exp(\pi \sqrt{-1}c_mx)\prod_{j=1}^{l_0+l_1} \theta(x+\tilde{t}_j(p))}{\theta(x)^{l_0}\theta(x+1/2)^{l_1}},
\label{tBVl0l1} \\
& \Lambda^{sym}_m(x,p)= \Lambda_m(x,p) -(-1)^{l_0}\Lambda_m(-x,p) .
\end{align}
By a similar argument to the case of the function $\Lambda^{sym}_{T}(x)$, the function $\Lambda^{sym}_m(x,p)$ is square-integrable on the interval $[0,1]$.
The function $\Lambda^{sym}_m(x,p)$ converges to the function $B'_m \Lambda^{sym}_{T}(x)$ $(c=l_0+l_1+2+2m)$ uniformly on compact sets for some non-zero constant $B'_m $ as $p\rightarrow 0$.

Combined with the relation $\Lambda^{sym}_{T}(x) = B_m \Phi(x)\psi _m(x)$ for some $B_m$, the following theorem is proved:
\begin{thm} \label{thm:pertu}
Let $m$ be an non-negative integer. For each $m$, there exists $\epsilon \in \Rea _{>0}$ such that if $|p|<\epsilon $ then there exists an eigenvalue $E_m(p)$ and an eigenfunction $g_m(x,p)$ of the Hamiltonian of the Inozemtsev model (see (\ref{Hamil})) with the condition $l_2=l_3=0$ such that 
\begin{align}
& E_m(p) \rightarrow \pi ^2 (l_0+l_1+2+2m)^2-\frac{\pi^2}{3}\left( l_0(l_0+1)+l_1(l_1+1) \right) \mbox{ and} \\
& g_m(x,p) \rightarrow  \Phi(x)\psi _m(x) \; \; \; \; \; \mbox{ on compact sets},
\end{align}
as $p\rightarrow 0$ and $g_m(x,p)$ is square-integrable, i.e.
\begin{equation}
\int_0^1 |g_m(x,p)|^2 dx<\infty.
\end{equation}
Here, the function $\psi _m(x)$ is defined in (\ref{Jacobipol}).
The functions $E_m(p)$ and $g_m(x,p)$ are holomorphic in $p$ if $|p|<\epsilon $.

\end{thm}
\begin{remk}
For the case $l_1=l_2=l_3=0$, this theorem was established in \cite{Tak}.
\end{remk}

\subsection{The algorithm of the perturbation}
Assume $l_2=l_3=0$. From formula (\ref{wpth}), we obtain the following expansion in $p$,
\begin{equation}
H=H_T-C_T+\sum_{k=1}^{\infty}V_k(x) p^{k},
\end{equation}
where $H$ is the Hamiltonian of the $BC_1$ (one-particle) Inozemtsev model (see (\ref{InoEF})) with the condition $l_2=l_3=0$, $H_T$ is the Hamiltonian of the $BC_1$ Calogero-Sutherland model (see (\ref{trigH})), $C_T$ is a constant, and $V_k(x)$ are the functions which are expressed as a finite sum of $\cos 2\pi n x$ $(n \in \Zint_{\geq 0})$.

Set $v_m=\Phi(x) \psi _m(x)$, then $v_m$ is an eigenvector of the operator $H_T$ with the eigenvalue $E_m=\pi ^2 (l_0+l_1+2+2m)^2$.

We now apply the algorithm of the perturbation.
More precisely, we determine the eigenvalues $\tilde{E}_m(p) = E_m+\sum_{k=1}^{\infty} \tilde{E}_{m}^{\{k\}}p^{k}$ and the eigenfunctions $v_m(p)= v_m+ \sum_{k=1}^{\infty} \sum_{m'} c_{m,m'}^{\{k\}}v_{m'}p^{k}$ of the operator $H_T -C_T +\sum_{k=1}^{\infty} V_k(x) p^{k}$ as formal power series in $p$.

Since the function $V_k(x)v_m$ is expressed as the finite sum of $v_{m'}$ $(m'\in \Zint_{\geq 0})$ and the eigenvalues $E_m$ $(m \in \Zint_{\geq 0})$ are non-degenerate, the coefficients of the formal power series $v_m(p)$ and $\tilde{E}_m(p)$ satisfying
\begin{align}
& \left( H_T -C_T +\sum_{k=1}^{\infty} V_k(x) p^{k} \right) v_m(p)= \tilde{E}_m(p)v_m(p), \label{Peqn} \\
& \langle v_m(p) , v_m(p)\rangle = \langle v_m , v_m\rangle ,\nonumber
\end{align}
are determined recursively and uniquely (for details see \cite{Tak2}).

The convergence of the formal power series $v_m(p)$ and $\tilde{E}_m(p)$ is not guaranteed a priori and they may diverge in general.
Nevertheless it is shown that the formal power series $v_m(p)$ and $\tilde{E}_m(p)$ converge in our model as a corollary of Theorem \ref{thm:pertu}.

\begin{cor}
Let $v_m(p)$ be the eigenvector and $\tilde{E}_m(p)$ the eigenvalue of the Hamiltonian $H$ of the $BC_1$ (one-particle) Inozemtsev model with the condition $l_2=l_3=0$, which are obtained by the algorithm of the perturbation.
If $|p|$ is sufficiently small then the power series $\tilde{E}_m(p)$ converges and the power series $v_m(p)$ converges on compact sets of $x$.
\end{cor}
\begin{proof}
From Theorem \ref{thm:pertu}, for each $m$ there exists $\epsilon \in \Rea _{>0}$ such that if $|p|<\epsilon $ then there exists the eigenvalue $E_m(p)$ and the square-integrable eigenfunction $g_m(x,p)$ that converge respectively to the trigonometric eigenvalue and eigenfunction as $p\rightarrow 0$.
Set $\tilde{g}_m(x,p)= \sqrt{\frac{\langle\Phi(x) \psi _m(x), \Phi(x) \psi _m(x)\rangle}{\langle g_m(x,p) , g_m(x,p)\rangle }}g_m(x,p)$, then $E_m(p)$ and $\tilde{g}_m(x,p)$ are holomorphic in $p$ near $p=0$ and they satisfy 
\begin{align}
& \left( H_T -C_T +\sum_{k=1}^{\infty} V_k(x) p^{k} \right) \tilde{g}_m(x,p)= E_m(p)\tilde{g}_m(x,p) ,\\
& \langle \tilde{g}_m(x,p) , \tilde{g}_m(x,p)\rangle = \langle\Phi(x) \psi _m(x), \Phi(x) \psi _m(x)\rangle = \langle v_m, v_m \rangle  . \nonumber
\end{align}
These equations are the same as (\ref{Peqn}). From the uniqueness of the coefficients in $p$, we obtain $\tilde{E}_m(p)= E_m(p)$ and $v_m(p)= \tilde{g}_m(x,p)$.
From the holomorphy of $E_m(p)$ and $\tilde{g}_m(x,p)$, the convergence of $\tilde{E}_m(p)$ and $v_m(p)$ in $p$ is shown.
\end{proof}

\section{Comments} \label{sec:com}
\subsection{}
After obtaining Theorem \ref{thm:pertu}, we should discuss the completeness of functions $g_m(x,p)$ (see Theorem \ref{thm:pertu}) and clarify properties of eigenvalues and eigenstates. These matters will be discussed in the forthcoming article \cite{Tak2}.
Note that some results will be obtained by applying the method performed in \cite{KT}.

\subsection{}
Roughly speaking, the problem of finding square-integrable eigenstates for the Hamiltonian of the Inozemtsev model (see (\ref{Ino})) is equivalent to the problem of finding the solutions to the Heun equation (see (\ref{eq:Heun})) with special conditions for a monodromy matrix around two singular points.

It would be desirable if the Bethe Ansatz method would be of help in investigating the monodromies of the Heun equation. 

\vspace{.3in}

{\bf Acknowledgment}
The author would like to thank Prof.~M. Kashiwara, Dr.~Y. Komori, and Prof.~T. Miwa for discussions and support. He is partially supported by the Grant-in-Aid for Scientific Research (No. 13740021) from the Japan Society for the Promotion of Science.

\appendix
\section{}
This appendix presents the definitions and formulae for elliptic functions.

Let $\omega_1$ and $\omega_3$ be complex numbers such that the value $\omega_3/ \omega_1$ is an element of the upper half plane.

The Weierstrass $\wp$-function, the Weierstrass sigma-function and the Weierstrass zeta-function are defined as follows,
\begin{align}
& \wp (z)=\wp(z|2\omega_1, 2\omega_3)= \\
& \; \; \; \; \frac{1}{z^2}+\sum_{(m,n)\in \Zint \times \Zint \setminus {(0,0)}} \left( \frac{1}{(z-2m\omega_1 -2n\omega_3)^2}-\frac{1}{(2m\omega_1 +2n\omega_3)^2}\right), \nonumber \\
& \sigma (z)=\wp(z|2\omega_1, 2\omega_3)=z\prod_{(m,n)\in \Zint \times \Zint \setminus \{(0,0)\} } \left(1-\frac{z}{2m\omega_1 +2n\omega_3}\right) \cdot \nonumber \\
& \; \; \; \;\cdot \exp\left(\frac{z}{2m\omega_1 +2n\omega_3}+\frac{z^2}{2(2m\omega_1 +2n\omega_3)^2}\right), \nonumber \\
& \zeta(z)=\frac{\sigma'(z)}{\sigma (z)}. \nonumber
\end{align}
Setting $\omega_2=-\omega_1-\omega_3$ and $e_i=\wp(\omega_i)$, $\eta_i=\zeta(\omega_i)$ $(i=1,2,3)$ yields the relations
\begin{align}
& e_1+e_2+e_3=\eta_1+\eta_2+\eta_3=0, \; \; \; \\
& \wp(z)=-\zeta'(x), \; \; \; (\wp'(z))^2=4(\wp(z)-e_1)(\wp(z)-e_2)(\wp(z)-e_3), \nonumber \\
& \wp(z+2\omega_j)=\wp(z), \; \; \; \zeta(z+2\omega_j)=\zeta(z)+2\eta_j ,\; \; \; \; (j=1,2,3), \nonumber \\
& \frac{\wp''(z)}{(\wp'(z))^2}=\frac{1}{2}\left( \frac{1}{z-e_1}+\frac{1}{z-e_2}+\frac{1}{z-e_3} \right), \nonumber \\
& \wp(z+\omega_i)=e_i+\frac{(e_i-e_{i'})(e_i-e_{i''})}{\wp(z)-e_i}, \; \; \; \; (i=1,2,3),\nonumber
\end{align}
where $i', i'' \in \{1,2,3\}$ with $i'<i''$, $i\neq i'$, and $i\neq i''$.
The co-sigma functions $\sigma_i(z)$ $(i=1,2,3)$ are defined by
\begin{align}
& \sigma_i(z)=\exp (-\eta_i z)\sigma(z+\omega_i)/\sigma(\omega _i),
\end{align}
then we obtain
\begin{align}
& \left( \frac{\sigma_i(z)}{\sigma(z)} \right )^2 =\wp(z)-e_i,\; \; \; \; \; \; \; (i=1,2,3). 
\end{align}

Let $\tau$ be an complex number on the upper half plane. The elliptic theta function $\theta (x)$ is defined by 
\begin{equation}
\theta(x)=2\sum_{n=1}^{\infty} (-1)^{n-1} \exp( \tau \pi \sqrt{-1}(n-(1/2))^2)\sin(2n-1)\pi x.
\end{equation}
Then
\begin{align}
& \theta(x+1)=-\theta (x), \; \; \; \; \theta(x+\tau)=-e^{-2\pi \sqrt{-1}x -\pi \sqrt{-1}\tau}\theta(x).
\end{align}
Setting $\tau=\omega_3/\omega_1$ yields the relation
\begin{equation}
\sigma(2\omega_1x)=2\omega_1 \exp(2\eta_1\omega_1x^2) \theta(x)/\theta '(0).
\label{sigthetaA}
\end{equation}
The following formulas are used to obtain Theorem \ref{thm:thetaBA}:
\begin{align}
& \frac{\theta '(\frac{1}{2})}{\theta(\frac{1}{2})}=0, \; \; \; \frac{\theta '(\pm \frac{ \tau}{2})}{\theta(\pm \frac{\tau}{2})}=\mp \pi \sqrt{-1},   \; \; \; \wp (x)=\frac{\theta '(x)^2}{\theta(x)^2}-\frac{\theta ''(x)}{\theta(x)}-2\eta_1,
\end{align}
where $\omega_1=1/2$ and $\omega_3=\tau /2$.

Setting $\omega_1=1/2$, $\omega_3=\tau /2$, and $p=\exp (2\pi \sqrt{-1} \tau )$, the expansion of the Weierstrass $\wp$ function in $p$ is given as follows:
\begin{equation}
\wp (x)=
\frac{\pi^2 }{\sin^2 (\pi x)}- \frac{\pi ^2}{3} -8\pi^2 \sum_{n=1}^{\infty} \frac{np ^{n}}{1-p ^{n}} (\cos 2n \pi x -1). \label{wpth}
\end{equation}

\end{document}